\documentclass[11pt,a4paper,reqno]{amsart}
\usepackage{amssymb,amsmath,amsthm}
\usepackage[utf8]{inputenc}
\usepackage[T1]{fontenc}
\usepackage{enumerate}
\usepackage[all]{xy}
\usepackage{fullpage}
\usepackage{comment}
\usepackage{array}
\usepackage{longtable}
\usepackage{stmaryrd}
\usepackage{mathrsfs} 
\usepackage{xcolor}
\usepackage{mathtools}
\usepackage{stmaryrd}
\SetSymbolFont{stmry}{bold}{U}{stmry}{m}{n}

\def\deg{\mathrm{deg}}

\def\int{\mathrm{Int}}

\usepackage{hyperref}
\hypersetup{
    colorlinks=true,
    linkcolor=blue,
    filecolor=red,  
    citecolor=green,
    urlcolor=cyan,
    pdftitle={Counterexamples to Minkowski's Conjecture and Escape of Mass in Positive Characteristic},
    pdfpagemode=FullScreen,
    }

\usepackage{cleveref}
\crefformat{section}{\S#2#1#3}
\crefformat{subsection}{\S#2#1#3}
\crefformat{subsubsection}{\S#2#1#3}
\usepackage{enumitem}
\usepackage{tikz}

\newtheorem{theorem}{Theorem}[section]
\newtheorem{lemma}[theorem]{Lemma}
\newtheorem{corollary}[theorem]{Corollary}
\newtheorem{proposition}[theorem]{Proposition}

\newtheorem{conjecture}[theorem]{Conjecture}

\theoremstyle{definition}
\newtheorem{definition}[theorem]{Definition}

\theoremstyle{remark}
\newtheorem{remark}[theorem]{Remark}

\makeatletter
\newcommand{\subalign}[1]{%
  \vcenter{%
    \Let@ \restore@math@cr \default@tag
    \baselineskip\fontdimen10 \scriptfont\tw@
    \advance\baselineskip\fontdimen12 \scriptfont\tw@
    \lineskip\thr@@\fontdimen8 \scriptfont\thr@@
    \lineskiplimit\lineskip
    \ialign{\hfil$\m@th\scriptstyle##$&$\m@th\scriptstyle{}##$\hfil\crcr
      #1\crcr
    }%
  }%
}
\makeatother

\numberwithin{equation}{section}

\title{Counterexamples to Minkowski's Conjecture and Escape of Mass in Positive Characteristic}
\author{Noy Soffer Aranov}
\email{noyso@campus.technion.ac.il}
\address{Department of Mathematics, Technion, Haifa, Israel}

\begin{document}
\begin{abstract}
    We show that there are infinitely many counterexamples to Minkowski's conjecture in positive characteristic regarding uniqueness of the upper bound of the multiplicative covering radius, $\mu$, by constructing a sequence of compact $A$-orbits where $\mu$ obtains its conjectured upper bound. In addition, we show that these orbits, as well as a slightly larger sequence of orbits, must exhibit complete escape of mass. 
\end{abstract}
\maketitle

\section{Introduction}
\label{sec:Intro}
Let $d\geq 2$ be an integer, let $G=\operatorname{SL}_d(\mathbb{R})$, let $\Gamma=\operatorname{SL}_d(\mathbb{Z})$, and let $X_d=G/\Gamma$. Then $X_d$ can be identified with the space of unimodular lattices in $\mathbb{R}^d$ through the identification $g\Gamma\mapsto g\mathbb{Z}^d$. Given a lattice $x\in X_d$ and a function $F:\mathbb{R}^d\rightarrow \mathbb{R}^+$, we define the $\operatorname{CovRad}_{F}(x)$ to be the infimal $r\geq 0$, such that for every $R>r$,
$$x+\{\mathbf{v}\in \mathbb{R}^d:F(\mathbf{v})<R\}=\mathbb{R}^d.$$
This value has been well studied for several functions $F$, such as the multiplicative function $N:\mathbb{R}^d\rightarrow \mathbb{R}^+$ defined by $N\left((v_1,\dots v_d)\right)=\prod_{i=1}^d\vert v_i\vert$. This function is dynamically significant, since it is invariant under the group of diagonal matrices with determinant $1$, which we denote by $A$. We define Minkowski's function as $\mu(x)=\operatorname{CovRad}_N(x)$. Since $N$ is $A$-invariant, then $\mu$ is $A$-invariant as well. Hence ergodicity of the $A$ action on $X_d$ implies that $\mu$ is constant almost everywhere, and in \cite{S09}, Shapira proved that for $d\geq 3$, $\mu(x)=0$ for Haar almost every $x\in X_d$. Furthermore, it is interesting to understand the set of values that $\mu$ obtains, and in particular, to understand the upper bound of $\mu$ is. A famous conjecture attributed to Minkowski claims the following:
\begin{conjecture}[Minkowski's Conjecture]
\label{realMink}
For every $d\geq 2$, and for every $x\in X_d$,
\begin{enumerate}
    \item \label{RMinkUpp} $\mu(x)\leq 2^{-d}=\mu(\mathbb{Z}^d)$.
    \item \label{RMinkUnique} $\mu(x)=2^{-d}$ if and only if $x\in A\mathbb{Z}^d$.
\end{enumerate}
\end{conjecture}
Conjecture \ref{realMink} has been proved for $d\leq 10$ (see for example \cite{Mink}, \cite{Dys}, \cite{Rem}, \cite{Sku}, \cite{HRS7}, \cite{HRS8}, \cite{KR9}, \cite{KR10}, and \cite{Sol}). Furthermore, in \cite{C}, Cassels proved that $2^{-d}$ is not isolated in the Minkowski spectrum
$$\mathcal{S}_d=\{\mu(x):x\in X_d\}$$
In fact Shapira proved a stronger fact, which relates to the structure of $A$-orbits. It is well known that $\mathbb{Z}^d$ has a divergent $A$-orbit, that is the function $a\in A\mapsto a\mathbb{Z}^d$ is a proper function. In particular $A\mathbb{Z}^d$ is not compact, but yet Cassels proved that $\mu(\mathbb{Z}^d)=2^{-d}$ can be approximated by evaluating $\mu$ at a sequence compact $A$-orbits.
\begin{theorem}[Main Theorem of \cite{C}]
\label{CasselsR}
There exists a sequence of compact $A$-orbits, $Ax_n\subseteq X_d$ such that $\mu(x_n)\rightarrow 2^{-d}$. 
\end{theorem}
The proof of Theorem \ref{CasselsR} is constructive, and it raises the following question - what can be the limit points of sequences of compact $A$-orbits? In \cite{S15}, Shapira provided a partial answer to this question by generalizing Cassels' construction. 
\begin{theorem}[Theorem 1.1 in \cite{S15}]
\label{mainS}
    For any $d\geq 2$, there exists a sequence of compact $A$-orbits $Ax_n\subseteq X_d$ such that any accumulation point of the form $x=\lim_{k\rightarrow \infty}a_kx_k$, where $a_k\in A$, must satisfy $x\in A\mathbb{Z}^d$.
\end{theorem}
Moreover, Shapira proved that the lattices satisfying the conclusion of Theorem \ref{mainS} must exhibit full escape of mass. It is well known that every compact $A$-orbit, $Ax_n$ supports a unique $A$-invariant probability measure $\mu_{Ax_n}$. We say that the $A$-orbits $Ax_n$ exhibit escape of mass if every limit point of $\mu_{Ax_n}$ gives mass $<1$ to $X_d$, and we say that the $A$-orbits $Ax_n$ exhibits full escape of mass if $\mu_{Ax_n}\rightarrow 0$.
\begin{corollary}[Corollary 1.2 in \cite{S15}]
    \label{EscMassS}
    The lattices satisfying the conclusion of Theorem \ref{mainS} must satisfy $\mu_{Ax_n}\rightarrow 0$.
\end{corollary}

In this paper, we shall prove a positive characteristic analogue of Theorem \ref{mainS}, as well as Corollary \ref{EscMassS}. This will lead to a positive characteristic analogue of Theorem \ref{CasselsR}, which will show that the conjectured upper bound of the Minkowski spectrum in positive characteristic is not unique, contrary to Conjecture \ref{realMink}(\ref{RMinkUnique}). 

\subsection{The Positive Characteristic Setting}
\label{subsec:Def}
We first introduce the positive characteristic setting. Let $d\geq 2$, let $p$ be a prime, $q$ be a power of $p$, and let $\mathcal{R}=\mathbb{F}_q[x]$ be the ring of polynomials over $\mathbb{F}_q$. Let $\mathcal{K}=\mathbb{F}_q(x)$ be the field of rational functions over $\mathbb{F}_q$. We define an absolute value on $\mathcal{R}$ by $\vert f\vert=q^{\deg(f)}$ and extend it to an absolute value on $\mathcal{K}$ by $\left|\frac{f}{g}\right|=q^{\deg(f)-\deg(g)}$. This absolute value satisfies the ultrametric inequality.
\begin{equation}
\label{eqn:ultrametric}
    \vert \alpha+\beta\vert\leq \max\{\vert \alpha\vert,\vert \beta\vert\}.
\end{equation}
Moreover, if $\vert \alpha\vert\neq \vert \beta\vert$, then, the inequality sign in (\ref{eqn:ultrametric}) is an equality. 

The topological completion of $\mathcal{K}$ with respect to the metric $d(f,g)=\vert f-g\vert$ is the field of Laurent series $\tilde{\mathcal{K}}$ defined by
$$\tilde{\mathcal{K}}\:= \mathbb{F}_q\left(\left(x^{-1}\right)\right)=\Bigg\{\sum_{n=-N}^{\infty} a_n x^{-n}:a_n\in \mathbb{F}_q,N\in \mathbb{Z}\Bigg\}.$$
Let $\mathcal{O}$ be the maximal compact subgroup of $\tilde{\mathcal{K}}$, that is
$$\mathcal{O}=\mathbb{F}_q\left[\left[x^{-1}\right]\right]=\{f\in \tilde{\mathcal{K}}:\vert f\vert\leq 1\}.$$
Denote by $\mathbf{U}$ the group of units, that is 
$$\mathbf{U}\:= \{f\in \tilde{\mathcal{K}}:\vert f\vert=1\}=\Bigg\{\sum_{n=0}^{\infty}a_nx^{-n}:a_n \in \mathbb{F}_q, a_0\in \mathbb{F}_q^*\Bigg\}=\mathcal{O}^*.$$
We can view $\tilde{\mathcal{K}}$ as the direct product $\tilde{\mathcal{K}}\cong \mathbb{Z}\times \mathbf{U}$ in the following way:
$$f\mapsto \left(\log_q\vert f\vert,\frac{f}{x^{\log_q\vert f\vert}}\right).$$
Define the functions $\rho(f)=\log_q\vert f\vert$ and $\pi(f)=\frac{f}{x^{\log_q\vert f\vert}}$. We often abuse notation and write $\rho(\mathbf{v})=(\rho(v_1)\dots \rho(v_d))$ and similarly $\pi(\mathbf{v})=(\pi(v_1),\dots \pi(v_d))$ for vectors $\mathbf{v}\in \tilde{\mathcal{K}}^d$. Similarly, for $g\in \operatorname{GL}_d(\tilde{\mathcal{K}})$ we define $(\rho(g))_{ij}=\rho(g_{ij})$ and $(\pi(g))_{ij}=\pi(g_{ij})$.

Let $G=\operatorname{GL}_d(\tilde{\mathcal{K}})$ be the group of invertible $d\times d$ matrices over $\tilde{\mathcal{K}}$ and let 
$$[G]=\operatorname{PGL}_d(\tilde{\mathcal{K}})\cong \operatorname{GL}_d(\tilde{\mathcal{K}})/\tilde{\mathcal{K}}^*I\cong \operatorname{GL}_d(\tilde{\mathcal{K}})/\tilde{\mathcal{K}}^*$$
be the group of invertible $d\times d$ matrices over $\tilde{\mathcal{K}}$ up to homothety. Let $\Pi: G\rightarrow [G]$ be the quotient map. Denote
$$[g]=\Pi(g)=\{cg:c\in \tilde{\mathcal{K}}^*\}.$$ 
Since $G$ is a topological group, $[G]$ inherits the quotient topology from $G$. Let $\Gamma=\operatorname{GL}_d(\mathcal{R})<G$ be the group of invertible $d\times d$ matrices with entries in $\mathcal{R}$ and let $[\Gamma]$ be its image under $\Pi$. Since $\Gamma$ is the stabilizer of $\mathcal{R}^d$ in $G$, then $[\Gamma]$ is the stabilizer of $\left[\mathcal{R}^d\right]$ in $[G]$ and thus, $[\Gamma]\cong \operatorname{GL}_d(\mathcal{R})/\tilde{\mathcal{K}}=\operatorname{PGL}_d(\mathcal{R})$. Let $\mathcal{L}_d=G/\Gamma$ and let $\left[\mathcal{L}_d\right]=[G]/[\Gamma]$. Since $[G]$ is a topological group and $[\Gamma]$ is a lattice in $[G]$ (see sections 2 and 3 of \cite{P}), then $\left[\mathcal{L}_d\right]$ inherits the quotient topology from $[G]$. Furthermore, $\left[\mathcal{L}_d\right]$ is identified with the space of lattices in $\tilde{\mathcal{K}}^d$ up to homothety via the identification
$$[g][\Gamma]\mapsto [g]\left[\mathcal{R}^d\right].$$
The determinant map $\det:G\rightarrow \tilde{\mathcal{K}}^*$ descends to a determinant map 
$$[\det]:[G]\rightarrow \tilde{\mathcal{K}}^*/(\tilde{\mathcal{K}}^*)^d.$$ 
through the quotient map $\Pi$. Since $\tilde{\mathcal{K}}\cong \mathbb{Z}\times \mathbf{U}$, then $$\tilde{\mathcal{K}}^*/(\tilde{\mathcal{K}}^*)^d\cong (\mathbb{Z}/d\mathbb{Z})\times (\mathbf{U}/\mathbf{U}^d).$$
Therefore, the image of $\big\vert [\det]\big\vert$ is $q^{\mathbb{Z}}/q^{d\mathbb{Z}}$. Thus, the set $\{1,q,q^2,\dots ,q^{d-1}\}$ is a set of representatives for $$\bigg\{\big\vert [\det]\left([g]\right)\big\vert:[g]\in [G]\bigg\}.$$ 
Let $[A]$ be the group of diagonal matrices in $[G]$ and let $A_1<G$ be the group of diagonal matrices $a$ with $\vert \det(a)\vert=1$. Let $[A_1]$ be the group of diagonal matrices $[a]\in [A]$ which have a representative $a'\in [a]$ with $\vert \det(a')\vert=1$. We identify $[A]$ with $A$, the group of matrices of determinants of absolute value lying in the set $\{1,q,q^2,\dots ,q^{d-1}\}$ by choosing a representative of every homothety class with the fitting determinant.

For $j=0,1,\dots, d-1$, we say that a lattice $\mathfrak{x}$ has determinant $q^j$ if there exists a representative of $\mathfrak{x}$ of the form $g\Gamma$ with $\vert \det(g)\vert=q^j$. We view $\mathcal{L}_d$ as $d$ copies of $\operatorname{SL}_d(\tilde{\mathcal{K}})$, every with determninant $q^j$ for $j=0,1,\dots,d$.
\begin{definition}
\label{length}
Given a lattice $\mathfrak{x}=[g][\Gamma]$, we define the length of the shortest non-zero vector in $\mathfrak{x}$ as
\begin{equation*}
    \ell(\mathfrak{x})=\frac{1}{\vert \det(g)\vert^{\frac{1}{d}}}\min\big\{\Vert \mathbf{v}\Vert:\mathbf{v}\in g\Gamma\setminus\{0\}\big\},
\end{equation*}
where $\Vert (v_1,\dots v_d)^t\Vert=\max_i\vert v_i\vert$. 
\end{definition}
In $\mathcal{L}_d$, Mahler's compactness criterion gives a necessary and sufficient condition for compactness (see \cite{Cas59} for the real case). Since $\left[\mathcal{L}_d\right]$ inherits the function $\ell:\left[\mathcal{L}_d\right]\rightarrow q^{\mathbb{Z}}$ from $\mathcal{L}_d$, then Mahler's compactness criterion also holds in $\left[\mathcal{L}_d\right]$ (see \cite{KST} for a version of Mahler's compactness criterion for general $S$-adic fields).
\begin{theorem}[Mahler's Compactness Criterion]
\label{Mahler}
A set of lattices $Y\subseteq \left[\mathcal{L}_d\right]$ is compact if and only if there exists $\varepsilon>0$ such that $\inf_{\mathfrak{x}\in Y}\ell(\mathfrak{x})>\varepsilon$. 

\end{theorem}
\begin{remark}
In the positive characteristic setting, we have to take lattices up to homothety instead of unimodular lattices, since there is no convenient normalization of lattices over $\tilde{\mathcal{K}}$. In $\mathbb{R}$, we can make any lattice $g\mathbb{Z}^d\subseteq \mathbb{R}^d$ unimodular by normalizing by $\vert\det(g)\vert^{1/d}$. On the other hand, if $\mathfrak{x}=[g][\Gamma]\in \left[\mathcal{L}_d\right]$ is a lattice, then $\det([g])$ may not necessarily have a $d$-th root in $\tilde{\mathcal{K}}$. For instance, if $\det([g])=x$, then $x^{1/d}\notin \tilde{\mathcal{K}}$. Therefore, it is more natural to work with lattices up to homothety. 
\end{remark}

\subsection{Main Results}
\label{subsec:mainThms}
Fix an integer $d\geq 2$ and a prime power $q$. We first define Minkowski's function in positive characteristic. Define the function $N:\tilde{\mathcal{K}}^d\rightarrow \mathbb{R}^+$ by
$$N(\mathbf{v})=\prod_{i=1}^d\vert v_i\vert.$$
We define $[\mathcal{G}_d]$ to be the space of translates of lattice, that is
$$[\mathcal{G}_d]=\big\{\mathfrak{x}+\mathbf{v}:\mathfrak{x}\in \left[\mathcal{L}_d\right],\mathbf{v}\in \tilde{\mathcal{K}}^d\big\}.$$
We identify $[\mathcal{G}_d]$ with the space 
$$\big\{g\mathcal{R}^d+\mathbf{v}:g\in G,\vert \det(g)\vert\in \{1,q,\dots ,q^{d-1}\},\mathbf{v}\in\tilde{\mathcal{K}}^d\big\}.$$
We define the projection $\operatorname{proj}:[\mathcal{G}_d]\rightarrow \left[\mathcal{L}_d\right]$ by $\mathfrak{x}+\mathbf{v}\mapsto \mathfrak{x}$. We identify the fiber $\operatorname{proj}^{-1}(\mathfrak{x})$ with the torus $\tilde{\mathcal{K}}^d/\mathfrak{x}$. Given $y=g\mathcal{R}^d+\mathbf{v}\in [\mathcal{G}_d]$, we define the product set of $y$ as
$$P(y)=\{N(\mathbf{w}):\mathbf{w}\in y\}=\big\{N(\mathbf{u}+\mathbf{v}):\mathbf{u}\in g\mathcal{R}^d\big\},$$
and we define $N(y)=\inf P(y)$. Given $\mathfrak{x}\in \left[\mathcal{L}_d\right]$, we define 
$$\mu(\mathfrak{x})=\operatorname{CovRad}_N(\mathfrak{x})=\frac{1}{\vert \det(g)\vert}\sup_{\mathbf{v}\in \tilde{\mathcal{K}}^d}\inf_{\mathbf{u}\in \mathfrak{x}}N(\mathbf{v}-\mathbf{u})=\frac{1}{\vert \det(g)\vert}\sup_{y\in \operatorname{proj}^{-1}(\mathfrak{x})}N(y),$$
where $g\mathcal{R}^d$ is a representative of the homothety class of $\mathfrak{x}$. We define the Minkowski spectrum by
$$\mathcal{S}_d=\{\mu(\mathfrak{x}):\mathfrak{x}\in [\mathcal{L}_d]\}.$$
It is easy to see that $\mu$ is $A$-invariant and $\mu\left(\left[\mathcal{R}^d\right]\right)=q^{-d}$. This enables us to conjecture what the upper bound of $\mu$ is.
\begin{conjecture}
    \label{MinkConj}
    For every $\mathfrak{x}\in \left[\mathcal{L}_d\right]$, $\mu(\mathfrak{x})\leq q^{-d}=\mu\left(\left[\mathcal{R}^d\right]\right)$.
\end{conjecture}
One natural question is whether the conjectured upper bound of $\mu$ is unique to the $[A]$-orbit of $\left[\mathcal{R}^d\right]$. That is, if $\mu(\mathfrak{x})=q^{-d}$, then, is it true that $\mathfrak{x}\in [A]\left[\mathcal{R}^d\right]$? In this paper we shall show that in contrast to the real case, $q^{-d}$, the conjectured upper bound of $\mu$, is not unique to $[A]\left[\mathcal{R}^d\right]$. Moreover, we prove a stronger claim.
\begin{theorem}
\label{Cassels}
There exist infinitely many compact $[A]$-orbits $[A]\mathfrak{x}$ such that $\mu(\mathfrak{x})=q^{-d}$.
\end{theorem}
In order to prove Theorem \ref{Cassels}, we shall prove a positive characteristic analogue of Theorem \ref{CasselsR} and use discreteness of the absolute value around non-zero points as well as the fact that the product sets $P(y)$ satisfy the following inheritance lemma (see \cite{S09} for the real analogue).
\begin{lemma}[Inheritance]
\label{inheritance}
If $y,y_0\in \mathcal{G}_d$ are such that $y_0\in \overline{Ay}$, then, $\overline{P(y_0)}\subseteq \overline{P(y)}$. 
\end{lemma}
\begin{remark}
\label{uppSemiCont}
A consequence of the Lemma \ref{inheritance} is the upper semicontinuity of $\mu$, that is if $\mathfrak{x}_n\rightarrow \mathfrak{x}$ in $\mathcal{L}_d$, then $\limsup\mu(\mathfrak{x}_n)\leq \mu(\mathfrak{x})$. Moreover, if $\mathfrak{x}_0\in \overline{A\mathfrak{x}}$, then $\mu(\mathfrak{x}_0)\geq \mu(\mathfrak{x})$. Ergodicity of the $A$-action on $\mathcal{L}_d$ with respect to the Haar measure implies that $\mu$ is constant almost everywhere. Furthermore, upper semicontinuity of $\mu$ implies that the generic value of $\mu$ is its minimal value.
\end{remark}
In order to prove Theorem \ref{Cassels}, we shall prove a positive characteristic analogue of Theorem \ref{mainS}.
\begin{theorem}
\label{main}
    Let $d\geq 2$. Then, there exists a sequence of lattices $\mathfrak{x}_k\in \left[\mathcal{L}_d\right]$ such that 
    \begin{enumerate}
        \item $[A]\mathfrak{x}_k$ is compact for every $k$, and
        \item Any limit point of the form $\mathfrak{x}=\lim_{k\rightarrow \infty}a_k\mathfrak{x}_k$ with $a_k\in [A]$ satisfies $\mathfrak{x}\in [A]\left[\mathcal{R}^d\right]$. 
    \end{enumerate}
\end{theorem}
From Theorem \ref{main}, we can obtain the following corollary, which pertains to escape of mass. This corollary can be viewed as an analogue of Corollary 1.2 in \cite{S15}. 
\begin{corollary}
    \label{EscMass}
    Let $[A]\mathfrak{x}_k$ be a sequence of compact orbits satisfying the conclusion of Theorem \ref{main}. Let $\mu_{[A]\mathfrak{x}_k}$ be the unique $[A]$-invariant probability measure supported on $[A]\mathfrak{x}_k$. Then $\mu_{[A]\mathfrak{x}_k}$ converge to the zero measure. 
\end{corollary}
\begin{proof}[Proof of Corollary \ref{EscMass}]
If $\mu$ is an accumulation point of $\mu_{[A]\mathfrak{x}_k}$, then by Theorem \ref{main}, $\mu$ must be supported on $[A]\left[\mathcal{R}^d\right]$. By Poincare recurrence, the only probability measure supported on $[A]\left[\mathcal{R}^d\right]$ is the $0$ measure, and thus, $\mu=0$.
\end{proof}
In order to prove Theorem \ref{main}, we shall provide precise bounds on the rate of convergence and the rate of escape of mass. For $\delta>0$ define the compact sets
$$\left[\mathcal{L}_d\right]^{\geq \delta}=\big\{\mathfrak{x}\in \left[\mathcal{L}_d\right]:\ell(\mathfrak{x})\geq \delta\big\}.$$
In \cref{sec:EscMass}, we show that a certain family of compact orbits $\mathcal{F}$ satisfies the following conditions:
\begin{enumerate}
    \item \label{eqn:muBnd}
    $\mu_{[A]\mathfrak{x}}\left(\left[\mathcal{L}_d\right]^{\geq \delta}\right)\ll f([A]\mathfrak{x})$ and
    \item \label{eqn:DisDivOrb}
    $\forall \mathfrak{y}\in [A]\mathfrak{x}\cap \left[\mathcal{L}_d\right]^{\geq \delta}$, $d\left(\mathfrak{y},[A]\left[\mathcal{R}^d\right]\right)\ll g([A]\mathfrak{x})$,
\end{enumerate}
where $f,g:\mathcal{F}\rightarrow \mathbb{R}$ are explicit functions satisfying $f([A]\mathfrak{x}),g([A]\mathfrak{x})\rightarrow 0$ as we vary $[A]\mathfrak{x}\in \mathcal{F}$. From (\ref{eqn:muBnd}) and Theorem \ref{Mahler}, we deduce that the orbits in $\mathcal{F}$ must satisfy the conclusion of Corollary \ref{EscMass}. Furthermore, the lattices satisfying (\ref{eqn:DisDivOrb}) must satisfy the conclusion of Theorem \ref{main}. All of these results are stated and proved in an effective manner as done in \cite{S15}. 
\begin{remark}
    All of our results can be generalized for global fields, but we will state them for $\tilde{\mathcal{K}}$, to avoid technicalities and to ease notations. 
\end{remark}

\subsection{Structure of this Article}
In \cref{sec:Mass}, we shall prove Theorem \ref{main}. To do so, in \cref{sec:EscMass}, we provide geometric definitions of $A$-orbits and prove that orbits with this geometry exhibit complete escape of mass. Then, in \cref{sec:ConstLattice}, we shall construct lattices satisfying the properties which we defined in \cref{sec:EscMass}. All of our proofs are completely analogous with those of \cite{S15} and \cite{C}. In \cref{subsec:Cass}, we shall show that a specific subsequence of the lattices constructed in \cref{sec:ConstLattice} satisfies the conclusion of Theorem \ref{Cassels}. 
\subsection{Acknowledgements}
I would like to thank Uri Shapira for introducing this problem to me and for carefully reading drafts of this paper. Without him, this article wouldn't be possible. I would also like to thank the anonymous referees for their effort in reviewing this article and their useful comments. This work has received funding from the European Research Council (ERC) under the European Union’s Horizon 2020 Research and Innovation Program, Grant agreement no. 754475.

\section{Escape of Mass}
\label{sec:Mass}
In this section we develop the necessary concepts that will allow us to establish the topological and distributional statements claimed above for the sequences of compact $A$-orbits we construct in \cref{sec:ConstLattice}. The majority of this section is completely identical to \cite{S15}, besides a few technical differences. 
\subsection{Simplex Sets}
\label{sec:Simplex}
We will first introduce the notation of simplex sets, which will be useful for the subsequent parts. Let $\Vert\cdot\Vert$ denote the supremum norm on $\tilde{\mathcal{K}}^d$, i.e. $\Vert \mathbf{v}\Vert=\max_i\vert v_i\vert$.
\begin{definition}
A simplex set $\Phi$ in $A_1$ is a set of $d$ matrices $\mathbf{t}_1,\dots ,\mathbf{t}_d\in A_1$, such that
\begin{enumerate}
    \item The group generated by $\Phi=\{\mathbf{t}_1,\dots ,\mathbf{t}_d\}$ is a lattice in $A_1$, and
    \item $\prod_{i=1}^d \mathbf{t}_i=I$.
\end{enumerate}
The associated lattice is $\Gamma_{\Phi}:=\langle \Phi\rangle$. Let $n=d-1$. Define
$$\mathbb{R}_0^d:=\Big\{(v_1,\dots ,v_d)\in \mathbb{R}^d:\sum_{i=1}^dv_i=0\Big\}\cong \mathbb{R}^n$$
and $\mathbb{Z}_0^d=\mathbb{R}_0^d\cap \mathbb{Z}^d$. For convenience we often write matrices in $A_1$ as vectors. We embed $A_1$ in $\tilde{\mathcal{K
}}^d$ by identifying $\operatorname{diag}\{a_1,\dots,a_d\}$ with $(a_1,\dots ,a_d)\in\tilde{\mathcal{K}}^d$. Under this identification, we obtain that $\rho(A_1)=\mathbb{Z}_0^d$. For $\mathbf{v}\in \mathbb{Z}_0^d$ denote $\lceil \mathbf{v}\rceil_{\mathbb{R}_0^d}=\max_i v_i$. For $\mathbf{a}\in A_1$, define $\lceil \mathbf{a}\rceil=q^{\lceil\rho(\mathbf{a})\rceil_{\mathbb{R}_0^d}}$. Define $$\xi_{\Phi}:=\max_{\mathbf{t}\in \Phi}\lceil \rho(\mathbf{t})\rceil_{\mathbb{R}_0^d}.$$
Let $S_{\Phi}:=\operatorname{hull}\{\rho(\Phi)\}$ be the convex hull of $\rho(\Phi)$ in $\mathbb{R}_0^d$ and let $\mathbf{S}_{\Phi}:=\rho^{-1}\left(\frac{n}{2}S_{\Phi}\cap \mathbb{Z}_0^d\right)\subseteq A_1$. Let $S_{\Phi}^o$ be the interior of $S_{\Phi}$ in $\mathbb{R}_0^d$ and let $\mathbf{S}_{\Phi}^o:=\rho^{-1}\left(\frac{n}{2}S_{\Phi}^o\cap \mathbb{Z}_0^d\right)$. 

Define $A_1(\mathbf{U}):=A_1\cap \operatorname{GL}_d(\mathbf{U})$. Let $\mathcal{P}_n$ be the group of permutations on $n$ elements. Given a simplex set $\Phi=\{\mathbf{t}_1,\dots ,\mathbf{t}_d\}$, define
\begin{equation}
\label{eqn:w}
    \mathbf{w}:=\frac{1}{d}\sum_{l=1}^d(l-1)\rho(\mathbf{t}_l)\in \mathbb{R}_0^d.
\end{equation}
For $\tau\in \mathcal{P}_n$, let $\mathbf{w}_{\tau}\in \mathbb{R}_0^d$ be vector obtained by permuting the coordinates of $\mathbf{w}$ by $\tau$. Let $$W_{\Phi}:=\{\mathbf{w}_{\tau}:\tau\in \mathcal{P}_n\}\subseteq \mathbb{R}_0^d.$$
\end{definition}
The following covering claim from \cite{S15} will be an essential part of our proofs.
\begin{proposition}[Proposition 3.8 in \cite{S15}]
\label{CoveringR}
Let $\Phi$ be a simplex set in $A_1$. Then, 
\begin{enumerate}
    \item \label{eqn:CovRadR}
    $\mathbb{R}_0^d=\frac{n}{2}S_{\Phi}+\rho(\Gamma_{\Phi})$ and
    \item \label{eqn:BndryPntsR}
    $\mathbb{R}_0^d\setminus \left(\frac{n}{2}S_{\Phi}^o+\rho(\Gamma_{\Phi})\right)\subseteq W_{\Phi}+\rho(\Gamma_{\Phi})$.
    \item \label{eqn:AlmBndryPntsR} There exists a universal constant $c>0$ such that for every $\gamma\in (0,1)$, 
    $$\mathbb{R}_0^d\setminus \left((1-\gamma)\frac{n}{2}S_{\Phi}+\rho(\Gamma_{\Phi})\right)\subseteq B_{c\gamma\xi_{\Phi}}(W_{\Phi})+\rho(\Gamma_{\Phi}),$$
    where
    $$B_{c\gamma\xi_{\Phi}}(W_{\Phi}):=\bigg\{\mathbf{v}\in \mathbb{Z}_0^d:\inf_{\mathbf{u}\in W_{\Phi}}\lceil \mathbf{v}-\mathbf{u}\rceil_{\mathbb{R}_0^d}\leq c\gamma\xi_{\Phi}\bigg\}.$$
\end{enumerate}
\end{proposition}

By intersecting Proposition \ref{CoveringR}(\ref{eqn:CovRadR}) and (\ref{eqn:AlmBndryPntsR}) with $\mathbb{Z}_0^d$ and then pulling these claims back with $\rho$, we obtain the following covering lemma in $A_1$.
\begin{lemma}
\label{Covering}
Let $\Phi$ be a simplex set in $A_1$. Then,
\begin{enumerate}
    \item \label{eqn:CovRad} 
    $\Gamma_{\Phi}\mathbf{S}_{\Phi}=A_1.$
    \item  \label{eqn:AlmBndaryPnts}
    For $0<\gamma<1$, define
    $$\mathbf{S}_{\Phi}^{(\gamma)}:=\rho^{-1}\left((1-\gamma)\frac{n}{2}S_{\Phi}\right)\cap \mathbb{Z}_0^d.$$
    Then, there exists a constant $c>0$ such that for any $0<\gamma<1$,
    $$A_1\setminus \left(\mathbf{S}_{\Phi}^{(\gamma)}\cdot \Gamma_{\Phi}A_1(\mathbf{U})\right)\subseteq \rho^{-1}\left(B_{c\gamma\xi_{\Phi}}(W_{\Phi})\right)\Gamma_{\Phi}A_1(\mathbf{U}).$$ 
\end{enumerate}
\end{lemma}
\begin{proof}
By definition, $\rho(\mathbf{a}\mathbf{b})=\rho(\mathbf{a})+\rho(\mathbf{b})$. Moreover, notice that 
        $$\rho^{-1}(\mathbf{0})=\big\{\mathbf{a}=\operatorname{diag}\{a_1\dots a_d\}\in A_1:\log_q\vert a_i\vert=1,\forall i\in\{1,\dots, d\}\big\}=A_1(\mathbf{U}).$$
     We can now use these facts to prove (\ref{eqn:CovRad}) and (\ref{eqn:AlmBndaryPnts}).
    \begin{enumerate}
        \item By Proposition \ref{CoveringR}(\ref{eqn:CovRadR}),
        \begin{equation*}
            A_1=\rho^{-1}(\mathbb{R}_0^d)=\rho^{-1}\left(\frac{n}{2}S_{\Phi}+\rho(\Gamma_{\Phi})\right)=\rho^{-1}\left(\frac{n}{2}S_\Phi\right)\Gamma_{\Phi}=\mathbf{S}_{\Phi}\Gamma_{\Phi}.
        \end{equation*}
        \item By Proposition \ref{CoveringR}(\ref{eqn:AlmBndryPntsR}), 
        \begin{equation}
        \begin{split}    
        \label{eqn:BallIncl}
        \rho^{-1}\left(\mathbb{R}_0^d\setminus\left((1-\gamma)\frac{n}{2}S_{\Phi}+\rho(\Gamma_{\Phi})\right)\right)\\
        \subseteq \rho^{-1}\left(B_{c\gamma\xi_{\Phi}}(W_{\Phi})+\rho(\Gamma_{\Phi})\right)\\
        =\rho^{-1}(B_{c\gamma\xi_{\Phi}}(W_{\Phi}))\Gamma_{\Phi}A_1(\mathbf{U}).
        \end{split}
        \end{equation}
        Notice that 
        \begin{equation}
        \label{eqn:A_1minusS_phi^gamma}
        \begin{split}
            \rho^{-1}\left(\mathbb{R}_0^d\setminus \left((1-\gamma)\frac{n}{2}S_{\Phi}+\rho(\Gamma_{\Phi})\right)\right)\\
            =A_1\setminus\rho^{-1}\left((1-\gamma)\frac{n}{2}S_{\Phi}\right)\Gamma_{\Phi}A_1(\mathbf{U})\\
            =A_1\setminus \mathbf{S}_{\phi}^{(\gamma)}\Gamma_{\Phi}A_1(\mathbf{U}).
        \end{split}
        \end{equation}
        Hence, by plugging (\ref{eqn:A_1minusS_phi^gamma}) into (\ref{eqn:BallIncl}), we obtain that 
        \begin{equation*}
            A_1\setminus\mathbf{S}_{\Phi}^{(\gamma)}\Gamma_{\Phi}A_1(\mathbf{U})\subseteq \rho^{-1}(B_{c\gamma\xi_{\Phi}}(W_{\Phi}))\Gamma_{\Phi}A_1(\mathbf{U}).
        \end{equation*}
    \end{enumerate}
\end{proof}
\subsection{Escape of Mass and Geometry of the Space of Lattices}
\label{sec:EscMass}
In this section, we shall connect between the covering lemmas obtained in \cref{sec:Simplex} and the structure of the $A_1$-orbit. This will provide conditions ensuring that a sequence of lattices to exhibits escape of mass. Let
$$\Omega=\big\{\mathfrak{x}\in \left[\mathcal{L}_d\right]:A_1\mathfrak{x}\text{ is compact}\big\}.$$
For $\mathfrak{x}\in \Omega$, we say that a simplex set $\Phi$ is a simplex set for $\mathfrak{x}$ if $\Gamma_{\Phi}=\langle \Phi\rangle\subseteq \operatorname{stab}_{A_1}(\mathfrak{x})$. For $\mathfrak{x}\in \Omega$, denote $\Delta_{\mathfrak{x}}:=\operatorname{stab}_{A_1}(\mathfrak{x})$, and we define $\vert \Delta_{\mathfrak{x}}\vert$ to be the determinant of the lattice $\rho(\Lambda_{\mathfrak{x}})\leq \mathbb{R}_0^d$.
 
We shall extract information about the structure of a compact $A_1$-orbit, $A_1\mathfrak{x}$, given that the length of the shortest vector of $\mathfrak{x}$ is very short and that $\mathfrak{x}$ has a simplex set of a nice form. We will need the fact that every $A_1$-orbit intersects a fixed compact set $\mathcal{L}_d^{\geq q^{-d}}$.
\begin{theorem}
\label{longVec}
There exists a universal constant $\delta_0>0$ such that for any $\mathfrak{x}\in \mathcal{L}_d$, $A\mathfrak{x}\cap \mathcal{L}_d^{\geq \delta_0}\neq \emptyset$. Furthermore, $\delta_0$ can be taken to be $\geq q^{-d}$.
\end{theorem}
\subsubsection{Proof of Theorem \ref{longVec}}
Our proof is very similar to Margulis' proof of the analogous result in $\mathbb{R}^d$, which can be found in the appendix of \cite{TW}. We shall include the proof for completeness. In order to prove Theorem \ref{longVec}, we shall need an analogue of Minkowski's Second Theorem.
\begin{theorem}[Equation (25) in \cite{Mah}]
\label{Mink2nd}
    Let $\lambda_i(\mathfrak{x})$ be the successive minima of $\mathfrak{x}=g\mathcal{R}^d\in \mathcal{L}_d$, that is
    $$\lambda_i(\mathfrak{x})=\min\{r>0:\text{there exist }i\text{ linearly independent vectors in }\mathfrak{x}\text{ of norm }\leq r\}.$$
    Then,
    $$\vert\det(g)\vert=\prod_{i=1}^d\lambda_i(\mathfrak{x}).$$
\end{theorem}
We shall now use Theorem \ref{Mink2nd} to prove the following analogue of Proposition A.1 in \cite{TW}.
\begin{proposition}
\label{finLengthen}
    For $r>0$, denote $B_r=B(0,r)$. Then, there exists a finite set $F\subseteq A_1$ such that for every $g\in G$ with $\vert \det(g)\vert\in \{1,q,\dots,q^{d-1}\}$, there exists $\mathbf{f}\in F$ such that
    \begin{equation*}
        \forall 0\neq \mathbf{w}\in g\mathcal{R}^d\cap B_{q^{-1}}, \Vert \mathbf{fw}\Vert\geq q\Vert \mathbf{w}\Vert.
    \end{equation*}
\end{proposition}
\begin{proof}
    By Theorem \ref{Mink2nd}, for every $g\in G$ with $\vert \det(g)\vert\in \{1,q,\dots ,q^{d-1}\}$, $\operatorname{span}\{g\mathcal{R}^d\cap B_{q^{-1}}\}$ is a proper subspace. Thus, it suffices to show that there exists a finite set $F\subseteq A_1$ such that for any proper subspace $V\subseteq \tilde{\mathcal{K}}^d$, there exists some $\mathbf{f}\in F$ such that for every $0\neq \mathbf{v}\in V$, $\Vert \mathbf{fv}\Vert\geq q\Vert \mathbf{v}\Vert$.

    Since $V$ is a proper subspace, there exists some $i$ such that $\operatorname{span}\{\mathbf{e}_i\}\cap V=\{0\}$. We want to choose $i$ wisely so that for every $\mathbf{v}\in V$, $\Vert \mathbf{b}_i\mathbf{v}\Vert\geq q\Vert \mathbf{v}\Vert$, where $(\mathbf{b}_i)_{jj}=\begin{cases}
        x^{-(d-1)}&i=j\\
        x&\text{else}
    \end{cases}\in A_1$. For $\mathbf{v}=(v_1,\dots, v_d)\in V$, define
    \begin{equation}
        \label{eqn:M_v}
        M_{\mathbf{v}}=\big\{l\in\{1,\dots, d\}:\left|v_l\right|=\Vert \mathbf{v}\Vert\big\}.
    \end{equation}
    If there exists some $l\in\{1,\dots ,d\}$, such that for every $0\neq \mathbf{v}\in V$, $M_{\mathbf{v}}\neq \{l\}$, then, for every $\mathbf{v}\in V$, there exists some $j_{\mathbf{v}}\neq l$, such that $\Vert \mathbf{v}\Vert=\left|v_{j_{\mathbf{v}}}\right|$. Hence,
    \begin{equation*}
        \Vert \mathbf{b}_l\mathbf{v}\Vert=q\left|v_{j_{\mathbf{v}}}\right|=q\Vert \mathbf{v}\Vert.
    \end{equation*}
    Thus, it suffices to show that there exists some $i$ such that $M_{\mathbf{v}}\neq \{i\}$ for every $0\neq \mathbf{v}\in V$. Assume on the contrary that for every $i\in\{1,\dots, d\}$, there exists some $\mathbf{v}^{(i)}\in V$ such that $M_{\mathbf{v}^{(i)}}=\{i\}$. Then for every $i\in\{1,\dots ,d\}$ and for every $j\neq i$, we have $\left|v^{(i)}_j\right|<\left|v^{(i)}_i\right|=\big\Vert \mathbf{v}^{(i)}\big\Vert$. Hence if $\sigma\neq I$ is some permutation in $\mathcal{P}_d$, then,
    \begin{equation*}
        \left|\prod_{i=1}^d v^{(i)}_{\sigma(i)}\right|<\left|\prod_{i=1}^d v^{(i)}_i\right|.
    \end{equation*}
    Thus, by the equality case of the ultrametric inequality (\ref{eqn:ultrametric}), the matrix whose columns are $\mathbf{v}^{(1)},\dots,\mathbf{v}^{(d)}$ has determinant of absolute value
    \begin{equation*}
        \left|\prod_{i=1}^d v^{(i)}_i+\sum_{I\neq \sigma\in \mathcal{P}_d}(-1)^{\operatorname{sgn}(\sigma)}\prod_{i=1}^dv^{(i)}_{\sigma(i)}\right|=\prod_{i=1}^d\left|v^{(i)}_i\right|=\prod_{i=1}^d\big\Vert \mathbf{v}^{(i)}\big\Vert\neq 0.
    \end{equation*}
    Therefore, $\mathbf{v}^{(1)},\dots,\mathbf{v}^{(d)}\in V$ are linearly independent, which contradicts the fact that $V$ is a proper subspace of $\tilde{\mathcal{K}}^d$. Thus, there exists some $i$ such that for every $\mathbf{v}\in V$, $M_{\mathbf{v}}\neq \{i\}$. Hence, the set $F=\{\mathbf{b}_1,\dots,\mathbf{b}_d\}\subseteq A_1$ satisfies the conditions of Proposition \ref{finLengthen}.
\end{proof}
\begin{proof}[Proof of Theorem \ref{longVec}]
    Let $F=\{\mathbf{b}_1,\dots,\mathbf{b}_d\}$ as defined in Proposition \ref{finLengthen}. Let $\mathfrak{x}\in \mathcal{L}_d$ and let $\mathbf{v}\in \mathfrak{x}$ satisfy $\Vert \mathbf{v}\Vert=\min_{0\neq \mathbf{w}\in \mathfrak{x}}\Vert \mathbf{w}\Vert$. We shall show that there exists some $\mathbf{a}\in A_1$ such that $\mathbf{a}\mathfrak{x}\cap B_{q^{-d}}=\{0\}$. The radius $q^{-d}$ was chosen since $q^{-d}=q^{-1}\min_{\mathbf{a}\in F}\min_{i=1,\dots ,d}\left|(\mathbf{a})_{ii}\right|$. Thus, for every $\mathbf{u}\notin B_{q^{-1}}$ and for every $\mathbf{a}\in F$, $\mathbf{au}\notin B_{q^{-d}}$. 

    If $\mathfrak{x}\cap B_{q^{-d}}=\{0\}$, then, $\mathfrak{x}\in A_1\mathfrak{x}\cap \mathcal{L}_d^{\geq q^{-d}}\neq \emptyset$. Now assume that $\mathfrak{x}\cap B_{q^{-d}}\neq \{0\}$. Since $\mathfrak{x}\cap B_{q^{-1}}$ spans a proper subspace of $\tilde{\mathcal{K}}^d$, Proposition \ref{finLengthen} implies that there exists some $\mathbf{a}_1\in F$ such that $\Vert \mathbf{a}_1\mathbf{u}\Vert\geq q\Vert \mathbf{u}\Vert$ for every $\mathbf{u}\in \mathfrak{x}\cap B_{q^{-1}}$. We shall now use Proposition \ref{finLengthen} to define a sequence $\mathbf{a}_k\in F$, satisfying $\Vert \mathbf{a}_k\mathbf{v}\Vert\geq q\Vert \mathbf{v}\Vert$, for every $\mathbf{v}\in \mathbf{a}_{k-1}\dots \mathbf{a}_1\mathfrak{x}\cap B_{q^{-1}}$. 
    
    Assume that we have already chosen $\mathbf{a}_1,\dots ,\mathbf{a}_k$ and denote $\tilde{\mathbf{a}}_k=\mathbf{a}_k\mathbf{a}_{k-1}\cdots \mathbf{a}_1$. In addition, assume that for every $j\leq k$, $B_{q^{-d}}\cap \tilde{\mathbf{a}}_j\mathfrak{x}\neq \{0\}$, since otherwise, we obtain that $A_1\mathfrak{x}\cap B_{q^{-d}}=\{0\}$, and then we can terminate this algorithm. Use proposition \ref{finLengthen} to choose some $\mathbf{a}_{k+1}\in F$ such that for every $\mathbf{v}\in B_{q^{-1}}\cap \tilde{\mathbf{a}}_k\mathfrak{x}$, $\Vert \mathbf{a}_{k+1}\mathbf{v}\Vert\geq q\Vert\mathbf{v}\Vert$. Let $\mathbf{v}\in B_{q^{-d}}\cap \tilde{\mathbf{a}}_{k+1}\mathfrak{x}$. Then, there exists $\mathbf{u}\in \tilde{\mathbf{a}}_k\mathfrak{x}$ such that $\mathbf{v}=\mathbf{a}_{k+1}\mathbf{u}$. Moreover, since
    \begin{equation*}
        q^{-n}\Vert \mathbf{u}\Vert\leq \Vert \mathbf{a}_{k+1}\mathbf{u}\Vert=\Vert \mathbf{v}\Vert\leq q^{-d},
    \end{equation*}
    then, $\mathbf{u}\in B_{q^{-1}}$. Hence, \begin{equation}
    \label{eqn:BallCont}
        B_{q^{-d}}\cap \tilde{\mathbf{a}}_{k+1}\mathfrak{x}\subseteq \mathbf{a}_{k+1}\left(B_{q^{-1}}\cap \tilde{\mathbf{a}}_k\mathfrak{x}\right).
    \end{equation}
    Let $\mathbf{v}$ be the shortest non-zero vector in $\tilde{\mathbf{a}}_{k+1}\mathfrak{x}$. If $\Vert \mathbf{v}\Vert\leq q^{-d}$, then (\ref{eqn:BallCont}) implies that there exists $\mathbf{u}\in B_{q^{-1}}\cap \tilde{\mathbf{a}}_k\mathfrak{x}$ such that
    \begin{equation*}
        \ell(\tilde{\mathbf{a}}_{k+1}\mathfrak{x})=\Vert \mathbf{v}\Vert=\Vert \mathbf{a}_{k+1}\mathbf{u}\Vert\geq q\Vert \mathbf{u}\Vert\geq q\ell(\tilde{\mathbf{a}}_k\mathfrak{x}).
    \end{equation*}
    Therefore, this process strictly increases the length of the shortest vector in $\tilde{\mathbf{a}}_k\mathfrak{x}$, so that for $k$ large enough, $\ell(\tilde{\mathbf{a}}_{k+1}\mathfrak{x})\geq q^{-d}$. Hence, $B_{q^{-d}}\cap \tilde{\mathbf{a}}_{k+1}\mathfrak{x}=\{0\}$ so that $A_1\mathfrak{x}\cap \mathcal{L}_d^{\geq q^{-d}}\neq \emptyset$. Hence, $A\mathfrak{x}\cap \mathcal{L}_d^{\geq q^{-d}}\neq \emptyset$. 
\end{proof}

\subsubsection{The Structure of the \texorpdfstring{$A$}{}-Orbit}
From now on, let $\mathfrak{x}\in \Omega$, let $\Phi$ be a simplex set for $\mathfrak{x}$, and let $\Gamma_{\Phi}$ be the associated lattice. We can now interpret the results of \cref{sec:Simplex} in terms of the structure of the $[A_1]$-orbit.
\begin{lemma}
\label{OrbStr}
Let $\Phi$ be a simplex set for $\mathfrak{x}$ and let $\Gamma_\Phi$ be the corresponding lattice. Then, 
$$A_1\mathfrak{x}=\big\{\mathbf{a}\mathfrak{x}:\mathbf{a}\in \mathbf{S}_{\Phi}\big\}.$$
\end{lemma}
\begin{proof}
By Lemma \ref{Covering}, every $\mathbf{a}\in A_1$ can be written as $\mathbf{a}'\mathbf{t}$ where $\mathbf{t}\in \Gamma_{\Phi}$ and $\mathbf{a}'\in \mathbf{S}_{\Phi}$. Thus,
$$\mathbf{a}\mathfrak{x}=\mathbf{a}'\mathbf{t}\mathfrak{x}=\mathbf{a}'\mathfrak{x}\in \big\{\mathbf{a}\mathfrak{x}:\mathbf{a}\in \mathbf{S}_{\Phi}\big\}.$$
\end{proof}
We shall now use Theorem \ref{longVec} and Lemma \ref{OrbStr} to bound the length of the shortest vector in $\mathfrak{x}$ with respect to its simplex set. 
\begin{lemma}
\label{ShortVec}
Let $\mathfrak{x}\in\Omega$ and let $\Phi=\{\mathbf{t}_1,\dots ,\mathbf{t}_d\}$ be a simplex set for $\mathfrak{x}$. Then, 
\begin{equation*}
    \ell(\mathfrak{x})\gg q^{-\frac{n}{2}\xi_{\Phi}}.
\end{equation*}
\end{lemma}
\begin{proof}
Let $\mathbf{v}$ be a shortest non-zero vector in $\mathfrak{x}$. Then by Lemma \ref{OrbStr} and Theorem \ref{longVec}, there exists $\mathbf{a}\in A_1$ with $\mathbf{a}\in \mathbf{S}_{\Phi}$ such that
\begin{equation}
\label{eqn:lBnddelta0}
    q^{-d}\leq \Vert \mathbf{a}\mathbf{v}\Vert\leq\lceil \mathbf{a}\rceil\cdot \ell(\mathfrak{x}).
\end{equation}
Write $\rho(\mathbf{a})=\frac{n}{2}\sum_{i=1}^d \alpha_i\rho(\mathbf{t}_i)$ where $\sum_{i=1}^d \alpha_i=1$. Then, 
\begin{equation}
\label{eqn:aNorm}
    \lceil \mathbf{a}\rceil=q^{\lceil \rho(\mathbf{a})\rceil_{\mathbb{R}_0^d}}\leq q^{\frac{n}{2}\max_{i}\vert \alpha_i\vert\cdot \lceil\rho(\mathbf{t}_i)\rceil_{\mathbb{R}_0^d}}\leq q^{\frac{n}{2}\xi_{\Phi}}.
\end{equation}
Thus, by plugging (\ref{eqn:aNorm}) into (\ref{eqn:lBnddelta0}) we obtain that $\ell(\mathfrak{x})\gg q^{-\frac{n}{2}\xi_{\Phi}}$.
\end{proof}
Motivated by Lemma \ref{ShortVec}, we make the following definition:
\begin{definition}
\label{Mtight}
Let $\mathfrak{x}\in \Omega$ and $M>1$ and let $\Phi$ be a simplex set for $\mathfrak{x}$. We say that $\Phi$ is $M$-tight if $\ell(\mathfrak{x})\leq Mq^{-\frac{n}{2}\xi_{\Phi}}$. Denote by $\Omega_M$ the set of lattices $\mathfrak{x}\in \Omega$ with an $M$-tight simplex set $\Phi$. 
\end{definition}
We shall now reinterpret Lemma \ref{Covering} in terms of the structure of the $A_1$-orbit. 
\begin{proposition}
\label{MeasCon}
Let $\mathfrak{x}\in \Omega_M$ be a lattice with an $M$-tight simplex set $\Phi$, let $c$ be the constant from Lemma \ref{Covering}(\ref{eqn:AlmBndaryPnts}), and let $\kappa\in \left(0,1\right)$. Define:
\begin{equation}
\label{eqn:delta}
    \delta=Mq^{-\frac{n}{2}\vert \Delta_{\mathfrak{x}}\vert^{\frac{\kappa}{n}}},
\end{equation}
\begin{equation}
\label{eqn:gamma}
    \gamma=\frac{\vert \Delta_{\mathfrak{x}}\vert^{\frac{\kappa}{n}}}{\xi_{\Phi}},
\end{equation}
\begin{equation}
    \label{eqn:r}
    r=c\vert\Delta_{\mathfrak{x}}\vert^{\frac{\kappa}{n}}.
\end{equation}
Define $$\mathbf{W}_{\Phi,\kappa}:=\rho^{-1}(B_r(W_{\Phi})\cap \mathbb{Z}_0^d).$$ 
Then,
\begin{enumerate}
    \item \label{eqn:DivOrbDelta}
    $\Big\{\mathbf{a}\mathfrak{x}:\mathbf{a}\in \mathbf{S}_{\Phi}^{(\gamma)}\cdot \Gamma_\Phi\Big\}\subseteq \left[\mathcal{L}_d\right]^{<\delta}$.
    \item \label{eqn:CompPart}
    $\Big\{\mathbf{a}\in A_1:\mathbf{a}\mathfrak{x}\in \left[\mathcal{L}_d\right]^{\geq \delta}\Big\}\subseteq \mathbf{W}_{\Phi,\kappa}\cdot\Gamma_{\Phi}$.
    \item \label{eqn:MeaCompPart}
    $\mu_{[A_1]\mathfrak{x}}\left(\left[\mathcal{L}_d\right]^{\geq \delta}\right)\ll \vert \Delta_{\mathfrak{x}}\vert^{-1+\kappa}$.
\end{enumerate}
\end{proposition}
We note that $\left|\Delta_{\mathfrak{x}}\right|\ll \xi_{\Phi}^n$, and therefore, $0<\gamma\ll \xi_{\Phi}^{-1+\kappa}<1$.
\begin{proof}[Proof of Proposition \ref{MeasCon}]
Let $\mathbf{v}$ be a shortest non-zero vector in $\mathfrak{x}$ and let $\mathbf{a}\in \mathbf{S}_{\Phi}^{(\gamma)}$. Then, 
\begin{equation*}
\begin{split}
    \Vert \mathbf{a}\mathbf{v}\Vert\leq \lceil \mathbf{a}\rceil \cdot\Vert \mathbf{v}\Vert\leq  Mq^{-\frac{n}{2}\xi_{\Phi}}q^{\frac{n}{2}(1-\gamma)\xi_{\Phi}}\\
    =Mq^{-\frac{n}{2}\gamma\xi_{\Phi}}=Mq^{-\frac{n}{2}\vert \Delta_{\mathfrak{x}}\vert^{\frac{\kappa}{n}}}=\delta,
    \end{split}
\end{equation*}
which proves (\ref{eqn:DivOrbDelta}). By (\ref{eqn:DivOrbDelta}),
\begin{equation}
\label{eqn:CompactPartofOrb}
    \{\mathbf{a}\in A_1:\mathbf{a}\mathfrak{x}\in \left[\mathcal{L}_d\right]^{\geq \delta}\}\subseteq A_1\setminus\left(\mathbf{S}_{\Phi}^{(\gamma)}\cdot \Gamma_{\Phi}\right).
\end{equation}
Thus, by Lemma \ref{Covering}(\ref{eqn:AlmBndaryPnts}) and (\ref{eqn:CompactPartofOrb}), 
\begin{equation*}
    \{\mathbf{a}\in A_1:\mathbf{a}\mathfrak{x}\in \left[\mathcal{L}_d\right]^{\geq \delta}\}\subseteq \rho^{-1}\left(B_{c\gamma\xi_{\Phi}}(W_{\Phi})\right)\Gamma_\Phi=\mathbf{W}_{\Phi,\kappa}\Gamma_{\Phi}.
\end{equation*}
Thus,
\begin{equation*}
    \mu_{A_1\mathfrak{x}}\left(\left[\mathcal{L}_d\right]^{\geq \delta}\right)\ll\frac{r^n}{\vert \Delta_{\mathfrak{x}}\vert}\ll \vert \Delta_{\mathfrak{x}}\vert^{-1+\kappa}.
\end{equation*}
\end{proof}
We shall make the following definition, which pertains to the structure of the $A_1$-orbit during the times $\mathbf{W}_{\Phi,\kappa}$.
\begin{definition}
\label{boundRetTimes}
Let $\varepsilon>0$, $M,J>0$. Denote by $\Omega_M(\varepsilon,J)$ the set of $\mathfrak{x}\in \Omega_M$ with an $M$-tight simplex set $\Phi$, such that for any $\mathbf{w}\in \mathbf{W}_{\Phi,\kappa}$, there exist $g\in G$, and $\mathbf{a},\mathbf{a}'\in A_1$ such that 
\begin{enumerate}
    \item \label{eqn:LattClose} $\mathbf{w}\mathfrak{x}=\mathbf{a}g\mathbf{a}'\mathcal{R}^d$,
    \item $\Vert g-\mathrm{Id}\Vert=\max_{i,j}\vert g_{ij}-\mathrm{Id}_{ij}\vert\leq Jq^{-\vert \Delta_{\mathfrak{x}}\vert^{\varepsilon}}$, and
    \item $\lceil\mathbf{a}\rceil\leq q^r$.
\end{enumerate}
\end{definition}
We shall now show that a sequence of lattices $\{\mathfrak{x}_k\}$ in $\Omega_M(\varepsilon,J)$ with $\vert \Delta_{\mathfrak{x}_k}\vert\rightarrow \infty$ must satisfy the conclusion of Theorem \ref{main}. 
\begin{proposition}
\label{longVecInDivOrb}
Fix $M,J>0$ and $\varepsilon>0$. Then there exists $\delta>0$ such that for any $\kappa <\min\{n\varepsilon,1\}$, for all but finitely many $\mathfrak{x}\in \Omega_M(\varepsilon,J)$ and for every $y\in A_1\mathfrak{x}\cap \left[\mathcal{L}_d\right]^{\geq \delta}$,
\begin{equation*}
    d(y,A_1\mathcal{R}^d)\leq Jq^{-\frac{1}{2}\vert \Delta_{\mathfrak{x}}\vert^{\varepsilon}}.
\end{equation*}
\end{proposition}
\begin{proof}
Let $\delta$ be as in (\ref{eqn:delta}), let $\mathfrak{x}\in \Omega_M(\varepsilon,J)$ and let $y\in A_1\mathfrak{x}\cap \left[\mathcal{L}_d\right]^{\geq \delta}$. Then by Proposition \ref{MeasCon}(\ref{eqn:CompPart}), there exists some $\mathbf{w}\in \mathbf{W}_{\Phi,\kappa}$ such that $y=\mathbf{w}\mathfrak{x}$. Since $\mathfrak{x}\in \Omega_M(\varepsilon,J)$, then, there exist some $\mathbf{a},\mathbf{a}'\in A_1$ and $g\in G$ satisfying $\lceil \mathbf{a}\rceil\leq q^r$ and $\Vert g-\mathrm{Id}\Vert\leq J\vert\Delta_{\mathfrak{x}}\vert^{-\varepsilon}$ such that
\begin{equation*}
    y=\mathbf{w}\mathfrak{x}=\mathbf{a}g\mathbf{a}'\mathcal{R}^d=\mathbf{a}g\mathbf{a}^{-1}(\mathbf{a}\mathbf{a}'\mathcal{R}^d).
\end{equation*}
Then,  
\begin{equation}
\label{eqn:y,AR^dDis}
\begin{split}
    d_{\left[\mathcal{L}_d\right]}(y,A_1\mathcal{R}^d)=d_{\left[\mathcal{L}_d\right]}\left(\mathbf{a}g\mathbf{a}^{-1}\left(\mathbf{a}\mathbf{a}'\mathcal{R}^d\right),A_1\mathcal{R}^d\right)\\
    \leq d_{G}(\mathrm{Id},\mathbf{a}g\mathbf{a}^{-1})\leq Jq^{dr}q^{-\vert \Delta_{\mathfrak{x}}\vert^{\varepsilon}}=Jq^{dc\vert \Delta_{\mathfrak{x}}\vert^{\frac{\kappa}{n}}-\vert \Delta_{\mathfrak{x}}\vert^{\varepsilon}}\\
    =Jq^{\vert \Delta_{\mathfrak{x}}\vert^{\varepsilon}\left(dc\vert \Delta_{\mathfrak{x}}\vert^{\frac{\kappa}{n}-\varepsilon}-1\right)}.
\end{split}
\end{equation}
If $\kappa<n\varepsilon$, then for $\vert \Delta_{\mathfrak{x}}\vert$ large enough, $dc\vert \Delta_{\mathfrak{x}}\vert^{\frac{\kappa}{n}-\varepsilon}<\frac{1}{2}$. Thus, 
\begin{equation*}
    d(y,A_1\mathcal{R}^d)\leq Jq^{-{\frac{1}{2}}\vert \Delta_{\mathfrak{x}}\vert^{\varepsilon}}.
\end{equation*}
\end{proof}

\subsubsection{Generating Simplex Sets and Visit Times to \texorpdfstring{$\left[\mathcal{L}_d\right]^{\geq \delta}$}{Lg}}
Given a simplex set $\Phi$ for a lattice $\mathfrak{x}\in \Omega$, it is desirable to determine whether $\Phi$ generates $\Delta_{\mathfrak{x}}$. In practice, it can be difficult to determine this. Therefore, in this section we shall provide conditions ensuring that $\langle \Phi\rangle=\Delta_{\mathfrak{x}}$ and we shall also show that under certain conditions the number of visit times to the compact part of the $A_1$-orbit is large. 

Given a lattice $\mathfrak{x}\in \Omega$, we say that $\mathbf{t}\in A_1$ is a visit time to $\left[\mathcal{L}_d\right]^{\geq \delta}$ if $\mathbf{t}\mathfrak{x}\in \left[\mathcal{L}_d\right]^{\geq \delta}$. We shall first distinguish between two distinct visit times. Unlike the previous definitions, which were all identical to the real setting, this definition differs from the analogous definition in \cite{S15}, since $\tilde{\mathcal{K}}$ is totally disconnected, in contrast to $\mathbb{R}$.
\label{subsec:RetTimes}
\begin{definition}
\label{def:EquivRetTimes}
Let $\mathbf{t}_1,\mathbf{t}_2\in A_1$ be two visit times to $\left[\mathcal{L}_d\right]^{\geq \delta}$. We say that $\mathbf{t}_1$ and $\mathbf{t}_2$ are equivalent visit times if 
\begin{equation*}
    \mathbf{t}_1\in B_2\mathbf{t}_2\mathfrak{x}\text{  and  }\mathbf{t}_2\in B_1\mathbf{t}_1\mathfrak{x},
\end{equation*}
where $B_i\subseteq A_1$ is the ball of maximal radius around $\mathrm{I}$ such that
\begin{equation*}
    B_i\mathbf{t}_i\mathfrak{x}\subseteq [\mathcal{L}_d]^{\geq \delta}.
\end{equation*}
\end{definition}
\begin{remark}
\label{DiffRetTimes}
If $\mathbf{t}_1,\mathbf{t}_2\in A_1$ are equivalent return times, then, there exist $\mathbf{a}_i\in B_i$ such that $\mathbf{a}_1\mathbf{t}_2\mathfrak{x}=\mathbf{t}_1\mathfrak{x}$ and $\mathbf{a}_2\mathbf{t}_1\mathfrak{x}=\mathbf{t}_2\mathfrak{x}$. Let $r_i$ be the radius of $B_i$. Thus, there exists $\Vert\mathbf{a}_2\Vert\leq q^{r_2}$ and $\Vert \mathbf{a}_1\Vert\leq q^{r_1}$, such that $\mathbf{a}_2\mathbf{t}_1\mathbf{t}_2^{-1},\mathbf{a}_1\mathbf{t}_2\mathbf{t}_1^{-1}\in \operatorname{stab}_{A_1}(\mathfrak{x})$. In particular, there exist $\mathbf{s}_i\in B_i\mathbf{t}_i$ for $i\in\{1,2\}$, such that $\mathbf{s}_1\mathbf{s}_2^{-1}\in \operatorname{stab}_{A_1}(\mathfrak{x})$. 
\end{remark}
We need some definitions regarding simplex sets. We define the standard simplex sets for $k\geq 1$ as
\begin{equation}
\label{eqn:stdSimp}
\Phi_*^k:=\{\mathbf{b}_1^{k},\dots,\mathbf{b}_d^{k}\}\text{, where } \mathbf{b}_l:=\begin{pmatrix}
x\\
\vdots \\
x^{-n}\\
\vdots \\
x
\end{pmatrix}\text{\} } l\text{-th coordinate}.
\end{equation}
and denote $\Gamma_*:=\Gamma_{\Phi_*}$, $\Delta_*:=\Delta_{\Phi_*}$ and $\Delta_*^k:=\Delta_{\Phi_*^k}$. By equation (3.4) in \cite{S09},
$$\mathbf{S}_k:=\mathbf{S}_{\Phi_*^k}=\rho^{-1}\left(\frac{n}{2}\operatorname{hull}\left(\rho\left(\Phi_*^{k}\right)\right)\right)=\Bigg\{\mathbf{a}\in A_1:\lceil \mathbf{a}\rceil \leq q^{\frac{n}{2}k}\Bigg\}.$$
\begin{definition}
\label{almStd}
Let $C\in \mathbb{N}$. We say that a simplex set $\Phi$ is $(k,C)$-standard if there exist $\mathbf{c}_i\in A_1$, $i\in\{1,\dots,d\}$, with $\Vert \mathbf{c}_i\Vert\leq q^{C}$ such that
\begin{equation*}
    \Phi=\big\{\mathbf{b}_1^{k}\mathbf{c}_1,\dots,\mathbf{b}_d^{k}\mathbf{c}_d\big\}.
\end{equation*}
We say that the associated lattice $\Gamma_{\Phi}=\langle \Phi\rangle$ is a $(k,C)$-standard lattice. 

We denote by $\Omega_M^{(C)}$ the set of lattices $\mathfrak{x}\in \Omega_M$, such that there exists $k$ such that $\mathfrak{x}$ has a $(k,C)$-standard simplex set $\Phi$ which is $M$-tight. 
\end{definition}
\begin{theorem}
\label{retTimes}
Fix $M>0,C\geq 0$, $0<\kappa<1$ and $0<\delta<1$. Let $\mathfrak{x}_k\in \Omega_M^{(C)}$ be such that there exists an $M$-tight simplex set $\Phi_k$ for $\mathfrak{x}_k$ which is $(k,C)$-standard. Let $\Gamma_k=\langle \Phi_k\rangle$, let $W_k=W_{\Phi_k}$, $\mathbf{W}_{k,\kappa}=\mathbf{W}_{\Phi_k,\kappa}=\rho^{-1}(\mathbb{Z}_0^d\cap B_r(W_k))$, where $r$ is as in (\ref{eqn:r}) and let
\begin{equation*}
    W_{k}'=\Big\{\mathbf{w}_{\tau}\in W_{k}:\exists\mathbf{t}\in \rho^{-1}(B_r(\mathbf{w}_{\tau}))\subseteq \mathbf{W}_{k,\kappa}:\mathbf{t}\mathfrak{x}_k\in \left[\mathcal{L}_d\right]^{\geq \delta}\Big\}\subseteq W_k.
\end{equation*}
Then, there exists some $k_0$ depending on $\kappa$ such that for every $k\geq k_0$,
\begin{enumerate}
    \item \label{eqn:GammaInd}
    $\left[\Delta_{\mathfrak{x}_k}:\Gamma_k\right]\leq n!$
    \item \label{eqn:W'=W}
    If $W_{k}'=W_{k}$, then  $\Delta_{\mathfrak{x}_k}A_1(\mathbf{U})=\Gamma_kA_1(\mathbf{U})$.
    \item \label{eqn:NumVisTime}
    If $W_{k}'=W_{k}$, then, for any $\delta_1\in (0,1)$, there are at least $n!$ distinct visits to $\left[\mathcal{L}_d\right]^{\geq \delta_1}$.
\end{enumerate}
\end{theorem}
\begin{proof}
    To save on notation, we denote $\mathfrak{x}_k=\mathfrak{x}$ and assume that $k$ is large enough so that the conclusion of Proposition \ref{MeasCon} holds. Let $\delta$ be as in (\ref{eqn:delta}), and consider $\mathbf{t}\in A_1$ such that $\mathbf{t}\mathfrak{x}\in \left[\mathcal{L}_d\right]^{\geq \delta}$. By Proposition \ref{MeasCon}(\ref{eqn:CompPart}) and the definition of $W_{k}'$,
    \begin{equation}
    \label{eqn:stabCont}
        \mathbf{t}\Delta_{\mathfrak{x}}\subseteq \rho^{-1}\left(B_r(W_k')\right)\Gamma_k,
    \end{equation}
    where $r\asymp \vert \Delta_{\mathfrak{x}}\vert^{\frac{\kappa}{n}}$. Equation (\ref{eqn:stabCont}) implies that
    \begin{equation}
    \label{eqn:rhoStabCont}
        \rho(\mathbf{t})+\rho(\Delta_{\mathfrak{x}})\subseteq B_r(W_k')+\rho(\Gamma_k).
    \end{equation}
    We first show that for every $\tau\in \mathcal{P}_n$, the coset $\rho(\mathbf{t})+\rho(\Delta_{\mathfrak{x}})$ can contain at most one point of $B_r(\mathbf{w}_{\tau})+\rho(\mathbf{v})$ for $\mathbf{v}\in \Gamma_k$. Statement (\ref{eqn:GammaInd}) will follow from this claim along with Lemma 4.11 in \cite{S15} when applied to $\rho(\Delta_{\mathfrak{x}})$ and $\rho(\Gamma_k)$. 
    
    Assume that there exist $\tau\in \mathcal{P}_n'$ and $\mathbf{v}\in \Gamma_k$ such that
    \begin{equation*}
        \left| \left(\rho(\mathbf{t})+\rho(\Delta_{\mathfrak{x}})\right)\cap \left(B_r(\mathbf{w}_{\tau})+\rho(\mathbf{v})\right)\right|\geq 2.
    \end{equation*}
    Then, there exist $\mathbf{u}_i\in \Delta_{\mathfrak{x}}$ and $\lceil\mathbf{s}_i\rceil_{\mathbb{R}_0^d}\leq r$ such that $\rho(\mathbf{t})+\rho(\mathbf{u}_i)=\mathbf{s}_i+\rho(\mathbf{v})$, for $i\in\{1,2\}$. Hence, $\rho(\mathbf{u}_2)-\rho(\mathbf{u}_1)=\mathbf{s}_2-\mathbf{s}_1\in \Delta_{\mathfrak{x}}$ has norm at most $dr$. On the other hand, the distance between the balls composing $B_r(W_k)+\rho(\Gamma_k)$ is greater than or equal to 
    \begin{equation}
    \begin{split}
    \label{eqn:DistBalls}
        \inf_{\mathbf{v}_i\in \Gamma_k,\sigma,\tau\in \mathcal{P}_n}d\left(B_r\mathbf{w}_{\tau})+\rho(\mathbf{v}_1),B_r(\mathbf{w}_{\sigma})+\rho(\mathbf{v}_2)\right)\\
        =\inf_{\mathbf{v}_i\in \Gamma_k,\sigma,\tau\in \mathcal{P}_n,\lceil \mathbf{s}_i\rceil_{\mathbb{R}_0^d}\leq r}\left(\mathbf{w}_{\tau}+\mathbf{s}_1+\rho(\mathbf{v}_1),\mathbf{w}_{\sigma}+\mathbf{s}_2+\rho(\mathbf{v}_2)\right)\\
        \gg\inf_{\mathbf{v}\in \Gamma_k,\lceil\mathbf{s}\rceil\leq rd}\lceil \rho(\mathbf{v})+\rho(\mathbf{s})\rceil\gg k+C-dr.
    \end{split}
    \end{equation}
    Since $\vert \Delta_{\mathfrak{x}}\vert\ll\xi_{\Phi_k}^n\ll (k+C)^n$, then, (\ref{eqn:DistBalls}) is $\gg \vert \Delta_{\mathfrak{x}}\vert^{\frac{1}{n}}\left(1-dc\vert\Delta_{\mathfrak{x}}\vert^{-\frac{1-\kappa}{n}}\right)$, which is larger than $dr=dc\vert \Delta_{\mathfrak{x}}\vert^{\frac{\kappa}{n}}$ as $\vert \Delta_{\mathfrak{x}}\vert\rightarrow \infty$, since $0<\kappa<1$. Thus, we obtain a contradiction to (\ref{eqn:rhoStabCont}).
    
    Let $\mathbf{u}\in \Delta_{\mathfrak{x}}$. Assume that $\mathbf{tu}\in \rho^{-1}\left(B_r(\mathbf{w}_{\tau})\right)\mathbf{v}$ where $\mathbf{v}\in \Gamma_k$ and $\tau\in \mathcal{P}_n$. Then, $\rho(\mathbf{t})+\rho(\mathbf{u})$ is the unique point of $\rho(\mathbf{t})+\rho(\Delta_{\mathfrak{x}})$ which is in $B_r(\mathbf{w}_{\tau})+\rho(\mathbf{v})$. Hence, there exists a unique $\mathbf{s}\in B_r(\mathbf{w}_{\tau})$ such that $\rho(\mathbf{t})+\rho(\mathbf{u})=\mathbf{s}+\rho(\mathbf{v})$. Thus, we can apply Lemma 4.11 in \cite{S15} to $\rho(\Delta_{\mathfrak{x}})$ and $\rho(\Gamma_k)$ to obtain that $[\rho(\Delta_{\mathfrak{x}}):\rho(\Gamma_{k})]\leq n!$. Hence, $[\Delta_{\mathfrak{x}}:\Gamma_{k}]\leq n!$, which proves (\ref{eqn:GammaInd}).
    
    We now prove (\ref{eqn:W'=W}). Let $\mathbf{t}$ be such that $\mathbf{t}\mathfrak{x}\in \left[\mathcal{L}_d\right]^{\geq \delta_1}$. By (\ref{eqn:stabCont}), we can assume that $\mathbf{t}\in \rho^{-1}\left(B_r(\mathbf{w}_{\tau})\right)$ for some $\mathbf{w}_{\tau}\in W_k'$. We shall show that $\rho(\mathbf{t})+\rho(\Delta_{\mathfrak{x}})\subseteq B_r(\mathbf{w}_{\tau})+\rho(\Gamma_k)$, which together with the fact that every ball composing $B_r(\mathbf{w}_{\tau})+\rho(\Gamma_k)$ contains at most one point of $\rho(\mathbf{t})+\rho(\Delta_{\mathfrak{x}})$, will imply that $\rho(\Gamma_k)=\rho(\Delta_{\mathfrak{x}})$. 
    
    Assume on the contrary that there exists some $\tau\neq \sigma\in \mathcal{P}_n$ such that
    \begin{equation*}
        \left(\rho(\mathbf{t})+\rho(\Delta_{\mathfrak{x}})\right)\cap \left(B_r(\mathbf{w}_{\sigma})+\rho(\Gamma_k)\right)\neq \emptyset.
    \end{equation*}
    Then, there exists some $\mathbf{v}\in \Delta_{\mathfrak{x}}$ and $\lceil \mathbf{s}_1\rceil_{\mathbb{R}_0^d}\leq r$ such that
    \begin{equation}
    \label{eqn:rho(t)sigma}
        \rho(\mathbf{t})=\mathbf{w}_{\sigma}+\mathbf{s}_1+\rho(\mathbf{v}).
    \end{equation}
    Hence, 
    \begin{equation*}
        \rho(\mathbf{t})-\mathbf{w}_{\sigma}-\mathbf{s}_1\in \rho(\Delta_{\mathfrak{x}}).
    \end{equation*}
    On the other hand, since $\mathbf{t}\in \rho^{-1}\left(B_r (\mathbf{w}_{\tau})\right)$, then, there exists some $\lceil\mathbf{s}_2\rceil_{\mathbb{R}_0^d}\leq r$ such that
    \begin{equation}
    \label{eqn:rho(t)tau}
        \rho(\mathbf{t})=\mathbf{w}_{\tau}+\mathbf{s}_2.
    \end{equation}
    Hence, if we denote $\mathbf{s}=\mathbf{s}_2-\mathbf{s}_1$, then, $\lceil \mathbf{s}\rceil_{\mathbb{R}_0^d}\leq nr$ and 
    \begin{equation}
    \label{eqn:shortVecStab}
        \mathbf{w}_{\tau}-\mathbf{w}_{\sigma}+\mathbf{s}=\rho(\mathbf{v})\in \rho(\Delta_{\mathfrak{x}}).
    \end{equation} 
    Define $\tau'\in \mathcal{P}_n$ by $\tau'^{-1}(j)=\sigma^{-1}(j)-\tau^{-1}(j)\mod d$. Then,
    \begin{equation*}
    \begin{split}
        \mathbf{w}_{\tau'}+\mathbf{w}_{\tau}-\mathbf{w}_{\sigma}=\\
        \frac{1}{d}\sum_{j=1}^d\left(\tau'^{-1}(j)-1+\tau^{-1}(j)-\sigma^{-1}(j)\right)\rho(\mathbf{b}_j^k\mathbf{c}_j)\\
        =-\frac{1}{d}\sum_{j=1}^d\rho(\mathbf{b}_j^k\mathbf{c}_j)=0\mod \rho(\Gamma_k).
    \end{split}
    \end{equation*}
    Hence, $\mathbf{w}_{\tau'}+\mathbf{w}_{\tau}-\mathbf{w}_{\sigma}\in \rho(\Gamma_k)$ so that $\mathbf{w}_{\tau'}+\mathbf{w}_{\tau}-\mathbf{w}_{\sigma}\notin W_k+\rho(\Gamma_k)$. Moreover, since $\vert \Delta_{\mathfrak{x}}\vert\ll (k+C)^n$ and for every $\mathbf{w}\in W_k$, $\mathbf{w}\asymp \frac{n}{2}(k+C)$, then, 
    \begin{equation}
    \begin{split}
    \label{eqn:distWGamma}
        d\left(\mathbf{w}_{\tau'}+\mathbf{w}_{\tau}-\mathbf{w}_{\sigma},W_{k}+\rho(\Gamma_k)\right)=d\left(0,W_k+\rho(\Gamma_k)\right)\\
        \gg \frac{n}{2}(k+C)\gg \frac{n}{2}\vert \Delta_{\mathfrak{x}}\vert^{\frac{1}{n}}\gg c\vert \Delta_{\mathfrak{x}}\vert^{\frac{\kappa}{n}}=r.
    \end{split}
    \end{equation}
    Hence, $\mathbf{w}_{\tau'}+\mathbf{w}_{\tau}-\mathbf{w}_{\sigma}\notin B_r(W_k)+\rho(\Gamma_k)$. Since $W_k=W_k'$, then $\mathbf{w}_{\tau'}\in W_{k}'$. Thus, there exists $\mathbf{t}'$ such that
    \begin{equation}
    \label{eqn:tau'Pt}
        \mathbf{t}'\in \rho^{-1}\left(B_r(\mathbf{w}_{\tau'})\right) \text{   and    }\mathbf{t}'\mathfrak{x}\in \left[\mathcal{L}_d\right]^{\geq \delta_1}.
    \end{equation}
    Therefore, there exists some $\lceil \mathbf{s}'\rceil_{\mathbb{R}_0^d}\leq r$ such that $\rho(\mathbf{t}')=\mathbf{w}_{\tau}'+\mathbf{s}'$.  Thus, (\ref{eqn:shortVecStab}) and (\ref{eqn:tau'Pt}) imply that 
    $$\mathbf{w}_{\tau'}+\mathbf{w}_{\tau}-\mathbf{w}_{\sigma}=\rho(\mathbf{t}')-\mathbf{s}'+\mathbf{s}-\rho(\mathbf{v})$$
    If we write $\tilde{\mathbf{s}}=\mathbf{s}-\mathbf{s}'$, then, $\lceil \tilde{\mathbf{s}}\rceil_{\mathbb{R}_0^d}\leq 2nr$ so that  (\ref{eqn:distWGamma}) implies that
    \begin{equation*}
    \begin{split}
        d\left(\rho(\mathbf{t}'),W_k+\rho(\Gamma_k)\right)=d\left(\mathbf{w}_{\tau'}+\mathbf{w}_{\tau}-\mathbf{w}_{\sigma}+\tilde{\mathbf{s}}+\mathbf{v},W_k+\rho(\Gamma_k)\right)\\
        \gg \frac{n}{2}\vert \Delta_{\mathfrak{x}}\vert^{\frac{1}{n}}-2nr\gg c\vert\Delta_{\mathfrak{x}}\vert^{\frac{\kappa}{n}}=r.
    \end{split}
    \end{equation*}
    This shows that for all but finitely many $\mathfrak{x}\in \Omega^{(C)}_M$, 
    $$\rho(\mathbf{t}')\notin B_r+W_k+\rho(\Gamma_k)= B_r(W_{k})+\rho(\Gamma_k).$$
    Thus, by Proposition \ref{MeasCon}(\ref{eqn:CompPart}), $\mathbf{t}'\mathfrak{x}\in \left[\mathcal{L}_d\right]^{<\delta_1}$, which contradicts (\ref{eqn:tau'Pt}). It follows that $\mathbf{t}\Delta_{\mathfrak{x}}\subseteq \rho^{-1}\left(B_r(\mathbf{w}_{\tau})\right)\Gamma_k$, which proves (\ref{eqn:W'=W}). 
    
    We now prove (\ref{eqn:NumVisTime}). Let $\mathbf{w}_{\tau}\in W_{k}''=W_{k}'$ and let $\mathbf{s}_{\tau}\in B_r$ be such that $\rho^{-1}\left(\mathbf{w}_{\tau}+\mathbf{s}_{\tau}\right)\mathfrak{x}\subseteq \left[\mathcal{L}_d\right]^{\geq \delta_1}$. Assume that there exist $\sigma\in \mathcal{P}_n$ and $\mathbf{s}_{\sigma}\in B_r$ such that for every $\mathbf{t}_{\tau}\in \rho^{-1}\left(\mathbf{w}_{\tau}+\mathbf{s}_{\tau}\right)$ and $\mathbf{t}_{\sigma}\in \rho^{-1}\left(\mathbf{w}_{\sigma}+\mathbf{s}_{\sigma}\right)$, $\mathbf{t}_{\sigma}$ and $\mathbf{t}_{\tau}$ are equivalent visit times. Then by Remark \ref{DiffRetTimes}, there exist $\mathbf{s}_1\in \rho^{-1}\left(B_r(\mathbf{w}_{\tau})\right)$ and $\mathbf{s}_2\in \rho^{-1}\left(B_r(\mathbf{w}_{\sigma})\right)$ such that $\mathbf{s}_1^{-1}\mathbf{s}_2\in \Delta_{\mathfrak{x}}$. Thus, there exist $\mathbf{s}_1',\mathbf{s}_2'\in B_r$ with $\rho(\mathbf{s}_1)=\mathbf{w}_{\tau}+\mathbf{s}_1'$ and $\rho(\mathbf{s}_2)=\mathbf{w}_{\sigma}+\mathbf{s}_2'$. Therefore, by taking $\mathbf{s}=\mathbf{s}_2'-\mathbf{s}_1'$, then $\lceil \mathbf{s}\rceil_{\mathbb{R}_0^d}\leq dr$ and also,
    $$\mathbf{w}_{\sigma}-\mathbf{w}_{\tau}+\mathbf{s}=\rho(\mathbf{s}_1^{-1}\mathbf{s}_2)=\rho(\mathbf{s}_2)-\rho(\mathbf{s}_1)\in \rho\left(\Delta_{\mathfrak{x}}\right).$$
   However this results in a contradiction as shown above after following (\ref{eqn:shortVecStab}).
\end{proof}
In conclusion, we obtain the following:
\begin{corollary}
\label{conclusion}
Fix $M,J,\varepsilon>0$ and some $C\geq 0$. Then,
\begin{enumerate}
    \item \label{eqn:SatEscMass}
    By Proposition \ref{MeasCon}, any sequence of distinct compact orbits $A_1\mathfrak{x}_k$ where $\mathfrak{x}_k\in \Omega_M$ must satisfy the conclusion of Corollary \ref{EscMass}.
    \item \label{eqn:SatLimPt}
    By Proposition \ref{longVecInDivOrb}, any sequence of distinct compact orbits $A_1\mathfrak{x}_k\in \Omega_M(\varepsilon,J)$ must satisfy the conclusion of Theorem \ref{main}.
    \item \label{eqn:SatRetTimes}
    If $\mathfrak{x}_k\in \Omega^{(C)}_M\cap \Omega_M(\varepsilon,J)$ is a sequence of lattices with $M$-tight simplex set $\Phi_k$ generating the lattice $\Gamma_k$, such that the $\mathbf{a}'$ from Definition \ref{boundRetTimes} are uniformly bounded, then, by Theorem \ref{retTimes}, for all sufficiently large $k$, $\Gamma_k$ generates $\Delta_{\mathfrak{x}_k}$ modulo $A_1(\mathbf{U})$. Moreover, for any $\delta\in (0,1)$, there are at least $n!$ distinct visit times to $\left[\mathcal{L}_d\right]^{\geq \delta}$.
\end{enumerate}
\end{corollary}
\begin{proof}
Parts (\ref{eqn:SatEscMass}) and (\ref{eqn:SatLimPt}) are direct consequences of Propositions \ref{MeasCon} and \ref{longVecInDivOrb} respectively. Part (\ref{eqn:SatRetTimes}) follows after observing that if $\mathfrak{x}_k\in \Omega_M^{(C)}\cap \Omega_M(\varepsilon,J)$ and the $\mathbf{a}'$ from Definition \ref{boundRetTimes} are uniformly bounded, then $W_k=W_k'$.
\end{proof}

\subsection{Construction of the Lattices}
\label{sec:ConstLattice}
\subsubsection{The Polynomials}
\label{subsec:poly}
We first define polynomials, which we use to construct the lattices exhibiting escape of mass. For every $\eta>0$, we define the following set of vectors in $\mathcal{R}^d$:
$$\mathcal{R}^d(\eta):=\Bigg\lbrace \mathbf{Q}\in \mathcal{R}^d:\forall i\neq j, \frac{\vert Q_i-Q_j\vert}{\Vert \mathbf{Q}\Vert }\geq \eta \Bigg\rbrace$$
Now for every $\mathbf{Q}\in \mathcal{R}^d(\eta)$, we can define a polynomial:
$$P_{\mathbf{Q}}(T):=\prod_{i=1}^{d}(T-Q_i)-1$$
We now prove a positive characteristic analogue of Lemma 5.1 in \cite{S15}. 
\begin{lemma}
\label{polyRoots}
Fix $\eta>0$. Then, for all but finitely many $\mathbf{Q}\in \mathcal{R}^d(\eta)$, the polynomial $P_{\mathbf{Q}}$ is irreducible over $\tilde{\mathcal{K}}$ and has $d$ distinct roots $\theta_j=\theta_j(\mathbf{Q})$ which all lie in $\tilde{\mathcal{K}}$. Moreover, 
\begin{equation*}
    \theta_j=Q_j+O_{\eta}\left(\Vert \mathbf{Q}\Vert^{-n}\right).
\end{equation*}
\end{lemma}
\begin{proof}
Let $\mathbf{Q}\in \mathcal{R}^d(\eta)$, let $L$ be the splitting field of $P_{\mathbf{Q}}$ over $\tilde{K}$, and let $\theta\in L$ be a root of $P_{\mathbf{Q}}$. Then,
\begin{equation}
\label{eqn:P_QDef}
    \prod_{i=1}^{d}(\theta-Q_i)=1.
\end{equation}
By (1.1) on page 3 in \cite{PR}, the absolute value on $\tilde{K}$ can be extended to a unique absolute value on $L$, which satisfies the ultrametric inequality. We abuse notation slightly, and denote this absolute value on $L$ by $\vert \cdot \vert$. Thus, (\ref{eqn:P_QDef}) implies that
\begin{equation}
    \label{eqn:polyNorm}
    \prod_{i=1}^d \vert \theta-Q_i\vert=1.
\end{equation}
If there exist $i\neq j$ such that $\vert \theta-Q_i\vert\leq 1$ and $\vert \theta-Q_j\vert\leq 1$, then by the ultrametric inequality, $\vert Q_i-Q_j\vert\leq 1$. This results in a contradiction for large enough $\Vert \mathbf{Q}\Vert$ since $\mathbf{Q}\in \mathcal{R}^d(\eta)$. Thus, there exists a unique $j$, such that $\vert \theta-Q_j\vert\leq 1$. Denote this unique $j$ by $j_{\theta}$. Then, if $j\neq j_{\theta}$, 
\begin{equation}
\label{eqn:theta-Q_j}
    \theta-Q_j=(\theta-Q_{j_{\theta}})+(Q_{j_{\theta}}-Q_j).
\end{equation}
Since $\vert\theta-Q_{j_{\theta}}\vert\leq 1$, then for large enough $\Vert \mathbf{Q}\Vert$, $\vert \theta-Q_j\vert=\vert Q_{j_{\theta}}-Q_j\vert$. Therefore, the definition of $\mathcal{R}^d(\eta)$ implies that
\begin{equation}
\label{eqn:jNotj_0}
    \vert \theta -Q_j\vert =\vert Q_j-Q_{j_{\theta}}\vert\asymp_{\eta} \Vert \mathbf{Q}\Vert.
\end{equation}
Thus, 
\begin{equation}
\label{eqn:theta-Q_j_0}
    \vert \theta-Q_{j_{\theta}}\vert =\prod_{i\neq j_{\theta}}\vert \theta-Q_i\vert^{-1} \asymp_{\eta} \Vert \mathbf{Q}\Vert^{-(d-1)}=\Vert\mathbf{Q}\Vert^{-n}.
\end{equation}
Let $\theta_j$ for $j=1,\dots,d$ be the roots of $P_{\mathbf{Q}}$. We shall now show that $\theta\mapsto j_{\theta}$ is an injective map from the set of roots of $P_{\mathbf{Q}}$ to $\lbrace 1,\dots,d\rbrace$. Assume that the map $\theta\mapsto j_{\theta}$ is not injective. Write $P_{\mathbf{Q}}(T)=\prod_{j=1}^{d}(T-\theta_j)$. Then there exist $j_1\neq j_2$ such that $l\:= j_{\theta_{j_1}}=j_{\theta_{j_2}}$. Therefore, 
\begin{equation*}
\begin{split}
    1=\left|\prod_{j=1}^d(Q_l-Q_j)-1\right|=\vert  P_{\mathbf{Q}}(Q_l)\vert\\
    =\prod_{j=1}^{d}\vert Q_l-\theta_j\vert\ll_{\eta} \Vert \mathbf{Q}\Vert^{-2n+(n-2)}=\Vert \mathbf{Q}\Vert^{-(d+1)}.
\end{split}
\end{equation*}
For large enough $\Vert \mathbf{Q}\Vert$, this results in a contradiction. Thus, we can order the $\theta_j$ so that
\begin{equation}
\label{eqn:samej}
\vert \theta_j-Q_j\vert\asymp_{\eta}\Vert \mathbf{Q}\Vert^{-n}\text{  and }
\end{equation}
\begin{equation}
\label{eqn:diffj}
\forall l\neq j, \vert \theta_j-Q_l\vert \asymp_{\eta} \Vert \mathbf{Q}\Vert. 
\end{equation}
In particular, this shows that $P_{\mathbf{Q}}$ has $d$ distinct roots, since
\begin{equation*}
    \vert \theta_i-\theta_j\vert=\vert\theta_i-Q_j+Q_j-\theta_j\vert= \vert \theta_i-Q_j\vert\asymp_{\eta}\Vert\mathbf{Q}\Vert.
\end{equation*}
We shall now show that $\theta_j\in \tilde{\mathcal{K}}$ for every $j\in\{1,\dots,d\}$. It is well known that roots in $L\setminus \tilde{\mathcal{K}}$ come in conjugate sets of size at least $2$ (see Chapter 1.14 and Chapter 2 of \cite{Neu}). Thus, if $\theta_j\notin \tilde{\mathcal{K}}$, then there exists an automorphism $\tau:L\rightarrow L$, which preserves $\tilde{\mathcal{K}}$, and some $i\neq j$ such that $\tau(\theta_j)=\theta_i$. It is well known that an automorphism of an extension of a local field equipped with an appropriate norm is an isometry (see Theorem 1.1 in \cite{CassLocalFields}). Since $\tau(\theta_j)=\theta_i$, then $\vert \theta_i\vert=\vert\theta_j\vert$ so that
$$\vert \theta_j-Q_i\vert=\vert \tau(\theta_j-Q_i)\vert=\vert \theta_i-Q_i\vert.$$ 
But, since $i\neq j$, 
$$\Vert \mathbf{Q}\Vert^{-n}\asymp_{\eta}\vert \theta_j-Q_j\vert=\vert \tau(\theta_j)-Q_j\vert=\vert \theta_i-Q_j\vert\asymp_{\eta} \Vert\mathbf{Q}\Vert,$$
which is a contradiction for $\Vert \mathbf{Q}\Vert$ large enough. Thus, $i=j$, so that of the $\theta_j$ must all lie in $\tilde{\mathcal{K}}$.

We now prove that $P_{\mathbf{Q}}$ is irreducible over $\mathcal{K}$. Since $\mathbb{F}_q$ is a field, then $\mathcal{R}=\mathbb{F}_q[x]$ is a unique factorization domain. Thus, by Gauss' lemma it suffices to prove that $P_{\mathbf{Q}}$ is irreducible over $\mathcal{R}$. If $P_{\mathbf{Q}}$ is reducible over $\mathcal{R}$, then there exists some proper subset $I\subset \lbrace 1,\dots,d\rbrace $, such that $F(T)=\prod_{j\in I}(T-\theta_j)$ is a polynomial over $\mathcal{R}$. Let $l\in I$. Then, $0\neq F(Q_l)\in \mathcal{R}$. On the other hand,
$$0\neq \vert F(Q_l)\vert=\prod_{j\in I}\vert Q_l-\theta_j\vert=O_{\eta}\left(\Vert \mathbf{Q} \Vert^{\vert I \vert-1-n}\right).$$
Therefore, $\vert F(Q_l)\vert<1$ for $\Vert \mathbf{Q}\Vert$ large enough, which is a contradiction to the assumption that $F(T)$ has coefficients in $\mathcal{R}$. Thus, for $\Vert \mathbf{Q}\Vert$ large enough $P_{\mathbf{Q}}$ must be irreducible over $\mathcal{R}$ and thus also be irreducible over $\mathcal{K}$. 
\end{proof}
\subsubsection{The Lattices}
\label{subsec:Lattice}
Fix $\eta>0$ and let $\mathbf{Q}\in \mathcal{R}^d(\eta)$ be such that $\Vert \mathbf{Q}\Vert$ is large enough so that $\mathbf{Q}$ satisfies the conclusion of Lemma \ref{polyRoots}. Let $\theta$ be a root of $P_{\mathbf{Q}}$ and let $\mathbb{F}_{\mathbf{Q}}=\mathcal{K}(\theta)$. Then by Lemma \ref{polyRoots}, $\mathcal{K}<\mathbb{F}_{\mathbf{Q}}\leq \tilde{\mathcal{K}}$ is an extension of degree $d$ over $\mathcal{K}$. Moreover, by (\ref{eqn:samej}) and (\ref{eqn:diffj}), we can order the embeddings $\sigma_1,\dots,\sigma_d:\mathbb{F}_{\mathbf{Q}}\rightarrow \tilde{\mathcal{K}}$ so that for every $j\in\{1,\dots,d\}$, $\theta_j=\sigma_j(\theta)$ satisfies $\theta_j=Q_j+O_{\eta}(\Vert \mathbf{Q}\Vert^{-n})$. Let 
$$\boldsymbol{\sigma}=\begin{pmatrix}
\sigma_1\\
\vdots\\
\sigma_d
\end{pmatrix}:\mathbb{F}_{\mathbf{Q}}\rightarrow \tilde{\mathcal{K}}^d,$$ 
and let
$$\Lambda_{\mathbf{Q}}:=\operatorname{span}_{\mathcal{R}}\{1,\theta,\dots,\theta^n\}.$$ 
Let
$\mathfrak{x}_{\mathbf{Q}}:=[\boldsymbol{\sigma}(\Lambda_{\mathbf{Q}})]\in [\mathcal{L}_d]$. In order to conclude the proof of Theorem \ref{main} and Corollary \ref{EscMass}, we shall show that the lattices $\mathfrak{x}_{\mathbf{Q}}$ satisfy the conditions of Corollary \ref{conclusion}. 
\begin{proposition}
\label{latticesWork}
For any $\eta>0$, there exist $M,J>0$, $C\geq 0$, and $\varepsilon>0$ such that for all but finitely many $\mathbf{Q}\in \mathcal{R}^d(\eta)$, $\mathfrak{x}_{\mathbf{Q}}\in \Omega_{M}(\varepsilon,J)\cap \Omega_M^{(C)}$. Moreover, for all but finitely many $\mathbf{Q}\in \mathcal{R}^d(\eta)$, there exists an $M$-tight simplex set $\Phi_{\mathbf{Q}}$, such that $\Phi_{\mathbf{Q}}$ generates $\Delta_{\mathfrak{x}}$ up to $A_1(\mathbf{U})$, and for every $\delta_1\in (0,1)$, $\left[\mathcal{L}_d\right]^{\geq \delta_1}$ contains at least $n!$ distinct visit times of $A_1\mathfrak{x}_{\mathbf{Q}}$.
\end{proposition}
\begin{proof}
Let $\Vert \mathbf{Q}\Vert$ be large enough so that the conclusion of Lemma \ref{polyRoots} holds. We shall first compute the determinant of $\boldsymbol{\sigma}(\Lambda_{\mathbf{Q}})$. Notice that $\boldsymbol{\sigma}(\Lambda_{\mathbf{Q}})$ is the lattice spanned over $\mathcal{R}$ by the columns of the Vandermonde matrix $\left(\theta_j^{i-1}\right)$. Thus,
\begin{equation}
\label{eqn:vdM}
\vert\det\left(\boldsymbol{\sigma}(\Lambda_{\mathbf{Q}})\right)\vert=\left|\det\begin{pmatrix}
1&\theta_1&\theta_1^2&\dots&\theta_1^n\\
1&\theta_2&\dots&\dots &\theta_2^n\\
\vdots&\vdots&\ddots &\dots &\vdots \\
1&\theta_d &\dots &\dots & \theta_d^n
\end{pmatrix}\right|=\prod_{i<j}\left|\theta_j-\theta_i\right|.
\end{equation}
Since $\theta_j=Q_j+O_{\eta}(\Vert \mathbf{Q}\Vert^{-n})$, then by (\ref{eqn:samej}) and (\ref{eqn:diffj}),
\begin{equation}
    \label{eqn:thetaDiff}
    \vert \theta_i-\theta_j\vert=\vert \theta_i-Q_j+Q_j-\theta_j\vert=\vert \theta_i-Q_j\vert\asymp_{\eta}\Vert \mathbf{Q}\Vert.
\end{equation}
Thus, (\ref{eqn:vdM}) is equal to 
\begin{equation}
\label{eqn:covolx_Q}
    \det\left(\boldsymbol{\sigma}(\Lambda_{\mathbf{Q}})\right)\asymp_{\eta}\prod_{i<j}\Vert \mathbf{Q}\Vert\asymp_{\eta}\Vert \mathbf{Q}\Vert^{\begin{pmatrix}
    d\\
    2
    \end{pmatrix}}.
\end{equation}
Since $\mathbf{1}=(1,\dots,1)^t\in \mathfrak{x}_{\mathbf{Q}}$ then, due to Definition \ref{length},
\begin{equation}
\label{eqn:vecLen}
\ell(\mathfrak{x}_{\mathbf{Q}})\leq \ell\left(\mathbf{1}\right)\ll_{\eta} \frac{1}{\Vert \mathbf{Q}\Vert^{\frac{1}{d}\begin{pmatrix}
d\\
2
\end{pmatrix}}}=\Vert \mathbf{Q}\Vert^{-\frac{n}{2}}.
\end{equation}
The link to the diagonal group stems from the following relationship - for any $\alpha,\beta\in \mathbb{F}_{\mathbf{Q}}$, we have 
\begin{equation*}
    \operatorname{diag}\left(\boldsymbol{\sigma}(\alpha)\right)\cdot\boldsymbol{\sigma}(\beta)=\begin{pmatrix}
\sigma_1(\alpha)\sigma_1(\beta)\\
\vdots\\
\sigma_d(\alpha)\sigma_d(\beta)
\end{pmatrix}=\begin{pmatrix}
\sigma_1(\alpha\beta)\\
\vdots\\
\sigma_d(\alpha\beta)
\end{pmatrix}=\boldsymbol{\sigma}(\alpha\beta).
\end{equation*}
Note that $\Lambda_{\mathbf{Q}}$ is the ring $\mathcal{R}[\theta]$. Denote $\omega_l=\theta-Q_l$ and note that $\omega_l\in \mathcal{R}[\theta]$. Furthermore, $\omega_l$ is a unit in $\mathcal{R}[\theta]=\Lambda_{\mathbf{Q}}$ since $\prod_{l=1}^d\vert \omega_l\vert=\prod_{l=1}^d\vert \theta-Q_l\vert=1$. Therefore, $\omega_l\Lambda_{\mathbf{Q}}=\Lambda_{\mathbf{Q}}$. Thus, if $\beta\in \Lambda_{\mathbf{Q}}$, then 
\begin{equation*}
    \operatorname{diag}(\boldsymbol{\sigma}(\omega_l))\boldsymbol{\sigma}(\beta)=\boldsymbol{\sigma}(\omega_l\beta)\in \boldsymbol{\sigma}(\Lambda_{\mathbf{Q}}).
\end{equation*}
Due to Dirichlet's units theorem (see chapter 3 in \cite{CF}), the group of units in $\mathcal{O}_{\mathbf{F}_{\mathbf{Q}}}$ is of rank $d-1$. Thus, $\operatorname{stab}_{[A]}(\mathfrak{x}_{\mathbf{Q}})$ is a lattice of rank $d-1$ in $[A]$. Hence, by the first isomorphism theorem, $[A]/\operatorname{stab}_{[A]}(\mathfrak{x}_{\mathbf{Q}})\equiv [A]\mathfrak{x}_{\mathbf{Q}}$. Therefore, $\mathfrak{x}_{\mathbf{Q}}$ has a compact $[A]$-orbit. 

We now show that $\mathfrak{x}_{\mathbf{Q}}\in \Omega_M^{(C)}$ for some constants $M,C>0$ (see Definition \ref{almStd}). Denote $\mathbf{t}_l=\operatorname{diag}\left(\boldsymbol{\sigma}(\omega_l)\right)$ and note that $\Phi_{\mathbf{Q}}=\{\mathbf{t}_1,\dots\mathbf{t}_d\}$ is a simplex set, since (\ref{eqn:samej}) and (\ref{eqn:diffj}) imply that the matrix with $\mathbf{t}_l$ in its columns has rank $n$. Then (\ref{eqn:samej}) and (\ref{eqn:diffj}) imply that
\begin{equation}
\label{eqn:disBound}
\xi_{\Phi_{\mathbf{Q}}}=\max_l\log_q\lceil\mathbf{t}_l\rceil=\log_q\Vert \mathbf{Q}\Vert+O_{\eta}(1).
\end{equation}
To ease on notations, we denote $\xi_{\mathbf{Q}}:=\xi_{\Phi_{\mathbf{Q}}}$ and $\Delta_{\mathbf{Q}}=\Delta_{\mathfrak{x}_{\mathbf{Q}}}$. Then, (\ref{eqn:samej}), (\ref{eqn:diffj}) and (\ref{eqn:disBound}) imply that there exists some $C=C_{\eta}\geq 0$ and some diagonal matrices $\lceil \mathbf{c}_l\rceil\leq q^C$, depending on $\mathbf{Q}$, such that for every $l\in\{1,\dots,d\}$, we have
\begin{equation}
\label{eqn:almStd}
   \mathbf{t}_l=\mathbf{b}_l^{\xi_{\mathbf{Q}}} \mathbf{c}_l.
\end{equation}
Thus, for large enough $\Vert \mathbf{Q}\Vert$, $\Phi_{\mathbf{Q}}$ is $(\xi_{\mathbf{Q}},C)$ standard. By combining (\ref{eqn:vecLen}) and (\ref{eqn:disBound}), we  obtain that 
\begin{equation}
   \label{eqn:shortVecLen} 
   \ell(\mathfrak{x}_{\mathbf{Q}})\ll_{\eta} \Vert \mathbf{Q}\Vert^{-\frac{n}{2}}\ll_{\eta} q^{-\frac{n}{2}\xi_{\mathbf{Q}}}.
\end{equation}
Thus, there exists some $M>0$ such that $\Phi$ is an $M$-tight simplex set for $\mathfrak{x}_{\mathbf{Q}}$. Hence, $\mathfrak{x}_{\mathbf{Q}}\in \Omega_M^{(C)}$ for large enough $\Vert \mathbf{Q}\Vert$. We now show that for all but finitely many $\mathbf{Q}$, $\mathfrak{x}_{\mathbf{Q}}\in \Omega_M(\varepsilon,J)$ for some $\varepsilon,J>0$ (see Definition \ref{boundRetTimes}). Take the change of basis given by
\begin{equation*}
    \{1,\theta,\dots, \theta^n\}\rightarrow \{1,\omega_1,\omega_1\omega_2,\dots, \omega_1\cdot \dots \cdot \omega_n\}.
\end{equation*}
Then, the columns of the following matrix are the images under $\boldsymbol{\sigma}$ of the basis \\
$\{1,\omega_1,\omega_1\omega_2,\dots, \omega_1\cdot \dots \cdot \omega_n\}$:
\begin{equation}
\label{eqn:M_Q}
    M_{\mathbf{Q}}=\begin{pmatrix}
1&\sigma_1(\omega_1)&\sigma_1(\omega_1\omega_2)&\dots &\sigma_1(\omega_1\dots \omega_n)\\
1&\sigma_2(\omega_1)&\sigma_2(\omega_1\omega_2)&\dots &\sigma_2(\omega_1\dots \omega_n)\\
\vdots&\dots &\ddots&\dots &\vdots \\
1&\sigma_d(\omega_1)&\sigma_d(\omega_1\omega_2)&\dots &\sigma_d(\omega_1\dots \omega_n)
\end{pmatrix}.
\end{equation}
Hence, $\mathfrak{x}_{\mathbf{Q}}=M_{\mathbf{Q}}\mathcal{R}^d$. Notice that (\ref{eqn:samej}) and (\ref{eqn:diffj}) imply that
\begin{equation}
\label{eqn:unitProdBND}
    \left| \sigma_i\left(\prod_{l=1}^j \omega_l\right)\right| \asymp_{\eta}\begin{cases}
\Vert \mathbf{Q}\Vert^{j-d}& i\leq j\\
\Vert \mathbf{Q}\Vert^j& i\geq j+1
\end{cases}.
\end{equation}
It is convenient to introduce some notation. Given matrices $g\in G$ and $g'\in \operatorname{GL}_d(\mathbb{R})$, we say that $g\ll_{\eta}^{av}g'$ if there exists some $c'>0$ depending only on $\eta$ such that for every $i,j\in\{1,\dots,d\}$, $\vert g_{i,j}\vert\leq c' g'_{i,j}$.

Fix $0<\kappa<1$. Let  $r$ be from (\ref{eqn:r}). Then, for every $\mathbf{w}\in \mathbf{W}_{\mathbf{Q},\kappa}=\rho^{-1}\left(B_r(W_{\mathbf{Q}})\right)$, there exist $\tau\in \mathcal{P}_n$ and $\mathbf{a}\in A_1$ with $\lceil \mathbf{a}\rceil\leq q^{r}$ such that 
\begin{equation}
\label{eqn:rho(w)Form}
    \rho(\mathbf{w})=\mathbf{w}_{\tau}+\rho(\mathbf{a})+O_d(1).
\end{equation}
We assume without loss of generality that $\tau=Id$, so that
\begin{equation}
\label{eqn:w_Id}
    \mathbf{w}_{Id}=\frac{1}{d}\sum_{l=1}^d(l-1)\rho(\mathbf{t}_l)=\frac{1}{d}\sum_{l=1}^d(l-1)\left(\xi_{\mathbf{Q}}\rho(\mathbf{b}_l)+\rho(\mathbf{c}_l)\right).
\end{equation}
By (\ref{eqn:disBound}), (\ref{eqn:almStd}) and (\ref{eqn:rho(w)Form}), 
\begin{equation}
\label{eqn:wBND}
\begin{split}
    \mathbf{a}^{-1}\mathbf{w}\ll_{\eta}^{av} \prod_{l=1}^d \operatorname{diag}(\Vert \mathbf{Q}\Vert^{l-1},\dots,\Vert \mathbf{Q}\Vert^{-n(l-1)},\dots, \Vert \mathbf{Q}\Vert^{l-1})^{\frac{1}{d}}\\
    \ll_{\eta}^{av} \operatorname{diag}\left(\Vert \mathbf{Q}\Vert^{\frac{n}{2}},\dots,\Vert \mathbf{Q}\Vert^{\frac{n}{2}-(i-1)},\dots,\Vert \mathbf{Q}\Vert^{-\frac{n}{2}}\right)\\
    =\Vert \mathbf{Q}\Vert^{\frac{n}{2}}\operatorname{diag}\left(1,\Vert \mathbf{Q}\Vert^{-1},\dots,\Vert \mathbf{Q}\Vert^{-n}\right).
\end{split}
\end{equation}
Since $M_{\mathbf{Q}}$ is obtained by changing the variables of $\boldsymbol{\sigma}(\Lambda_{\mathbf{Q}})$, then, by (\ref{eqn:covolx_Q}), $\vert \det(M_{\mathbf{Q}})\vert^{\frac{1}{d}}\asymp_{\eta}\Vert \mathbf{Q}\Vert^{\frac{n}{2}}$. Let $\overline{M}_{\mathbf{Q}}$ be the representative in the homothety class of $M_{\mathbf{Q}}$ satisfying $\vert \det(\overline{M}_{\mathbf{Q}})\vert\in \{1,q,\dots,q^{n}\}$. Then (\ref{eqn:unitProdBND}) implies that
\begin{equation}
\label{eqn:y_QBND}
    \overline{M}_{\mathbf{Q}}\ll_{\eta}^{av} \Vert \mathbf{Q}\Vert^{-\frac{n}{2}}\begin{pmatrix}
1& \Vert \mathbf{Q}\Vert^{-n}& \Vert \mathbf{Q}\Vert^{-(n-1)}&\dots &\Vert \mathbf{Q}\Vert^{-1}\\
1&\Vert \mathbf{Q}\Vert&\Vert \mathbf{Q}\Vert^{-(n-1)}&\dots &\Vert \mathbf{Q}\Vert^{-1}\\
\vdots&\dots&\ddots &\dots &\vdots\\
1&\Vert \mathbf{Q}\Vert&\Vert \mathbf{Q}\Vert^2&\dots &\Vert \mathbf{Q}\Vert^n
\end{pmatrix}.
\end{equation}
Thus, by combining (\ref{eqn:wBND}) and (\ref{eqn:y_QBND}), we obtain that
\begin{equation}
    \mathbf{a}^{-1}\mathbf{w}\overline{M}_{\mathbf{Q}}\ll_{\eta}^{av}.
\begin{pmatrix}
1&\Vert \mathbf{Q}\Vert^{-n}&\Vert \mathbf{Q}\Vert^{-(n-1)}&\dots&\Vert \mathbf{Q}\Vert^{-1}\\
\Vert \mathbf{Q}\Vert^{-1}&1&\Vert \mathbf{Q}\Vert^{-d}&\dots&\Vert \mathbf{Q}\Vert^{-2}\\
\vdots&\dots&\ddots&\dots&\vdots\\
\Vert \mathbf{Q}\Vert^{-n}&\Vert \mathbf{Q}\Vert^{-(n-1)}&\dots&\dots&1
\end{pmatrix}.
\end{equation}
In conclusion, if we denote $(\mathbf{a}^{-1}\mathbf{w}\overline{M}_{\mathbf{Q}})_{ij}=c_{ij}$, then,
\begin{enumerate}
    \item \label{eqn:bndc}$\vert c_{ii}\vert \ll_{\eta}1$ for every $i\in\{1,\dots,d\}$,
    \item $\vert c_{ij}\vert\ll_{\eta}\Vert \mathbf{Q}\Vert^{-1}$ for $i\neq j$, and
    \item $1=\det(c_{ij})=\prod_{i=1}^d c_{ii}+O_{\eta}(\Vert \mathbf{Q}\Vert^{-1})$.
\end{enumerate}
Therefore, we can write
\begin{equation}
\label{eqn:decom}
(c_{ij})=\begin{pmatrix}
\frac{c_{11}}{c_{11}}&\frac{c_{12}}{c_{22}}&\dots &\frac{c_{1d}}{c_{dd}}\\
\frac{c_{21}}{c_{11}}&\frac{c_{22}}{c_{22}}&\dots &\frac{c_{2d}}{c_{dd}}\\
\vdots&\dots&\ddots&\dots\\
\frac{c_{d1}}{c_{11}}&\dots&\dots&\frac{c_{dd}}{c_{dd}}
\end{pmatrix}\cdot
\begin{pmatrix}
c_{11}& & & \\
& \ddots & &\\
& & \ddots &\\
& & & c_{dd}
\end{pmatrix}.
\end{equation}
Denote the right hand side of (\ref{eqn:decom}) as $g\cdot \mathbf{a}'$. Notice that
\begin{equation*}
    g_{ij}\ll_{\eta}\begin{cases}
1&i=j\\
\Vert \mathbf{Q}\Vert^{-1}&i\neq j
\end{cases}.
\end{equation*}
Therefore, $\Vert g-Id\Vert\ll_{\eta} \Vert \mathbf{Q}\Vert^{-1}\ll_{\eta} q^{-\xi_{\mathbf{Q}}}$ by (\ref{eqn:disBound}). Since $\vert \Delta_{\mathbf{Q}}\vert \ll_{\eta}\xi_{\mathbf{Q}}^n$, then by choosing $\varepsilon$ and $J$ appropriately, we can ensure that $\Vert g-Id\Vert\leq Jq^{-\vert \Delta_{\mathbf{Q}}\vert^{\varepsilon}}$. This completes the verification of Definition \ref{boundRetTimes} and thus, $\mathfrak{x}_{\mathbf{Q}}\in \Omega_M(\varepsilon,J)$. 

Finally, the remaining part of the statement follows from Corollary \ref{conclusion}(\ref{eqn:SatRetTimes}) after observing that the matrices $\mathbf{a}'$ from Definition \ref{boundRetTimes} are uniformly bounded.
\end{proof}
\section{Proof of Theorem \ref{Cassels}}
\label{subsec:Cass}
In this section, we shall show that a subsequence of the lattices we constructed in \cref{sec:ConstLattice} satisfy the conclusion of Theorem \ref{Cassels}. To do so, we shall prove that this subsequence $\{\mathfrak{x}_\mathbf{Q}\}_{\mathbf{Q}\in I}$ satisfies $\lim_{\mathbf{Q}\rightarrow\infty,\mathbf{Q}\in I}\mu(\mathfrak{x}_\mathbf{Q})=q^{-d}$.

We first return to the notations of \cref{sec:ConstLattice} to better understand the lattices $\mathfrak{x}_{\mathbf{Q}}$. Fix distinct polynomials $a_1,\dots,a_d$ such that $a_j\equiv 0 \mod x^2$. For $Q\in \mathcal{R}$, let $\mathbf{Q}:=(Qa_1,\dots,Qa_d)$. Assume that $\eta$ is small enough so that for every $Q\in \mathcal{R}$, $\mathbf{Q}\in \mathcal{R}^d(\eta)$. Let $\theta$ be a root of 
$$P_{\mathbf{Q}}(t)=\prod_{i=1}^d(t-Qa_i)-1.$$ 
Let $\mathfrak{x}_{\mathbf{Q}}=[\boldsymbol{\sigma}(\Lambda_{\mathbf{Q}})]$. Notice that (\ref{eqn:jNotj_0}) and (\ref{eqn:theta-Q_j_0}) imply that in the symbols of \cref{sec:ConstLattice}, 
\begin{equation}
    \label{eqn:diffjCass}
    \forall i\neq j, \vert \theta_i-Qa_j\vert\asymp_{\eta}\Vert \mathbf{Q}\Vert
\end{equation}
and
\begin{equation}
    \label{eqn:samejCass}
    \vert \theta_i-Qa_i\vert\asymp_{\eta}\Vert \mathbf{Q}\Vert^{-n}.
\end{equation}
Thus, (\ref{eqn:diffjCass}) and (\ref{eqn:samejCass}) along with (\ref{eqn:vdM}) imply together that 
\begin{equation*}
    \det(\mathfrak{x}_{\mathbf{Q}})\asymp_{\eta}\Vert \mathbf{Q}\Vert^{\begin{pmatrix}
    d\\
    2
    \end{pmatrix}}.
\end{equation*}
Let $\theta=\theta_{\mathbf{Q}}$ be a root of $P_{\mathbf{Q}}$ and let $\omega_i=\theta-Qa_i$. Since $P_{\mathbf{Q}}(\theta)=0$, then $\prod_{i=1}^d \omega_i=1$, and therefore, $\omega_i$ are units in $\mathcal{O}_{\mathbb{F}_{\mathbf{Q}}}$. Denote $u_1=1$ and for $j\in\{2,\dots,n\}$, denote $u_j=\omega_1\cdots \omega_{j-1}$. Since $\theta^k$ is a linear combination of $u_1,\dots,u_{k+1}$ with coefficients in $\mathcal{R}$, then $u_1,\dots, u_d$ is an $\mathcal{R}$ basis for $\Lambda_{\mathbf{Q}}=\mathcal{R}[\theta]$.  Notice that by (\ref{eqn:unitProdBND}), (\ref{eqn:samejCass}) and (\ref{eqn:diffjCass}), for every $j\in\{1,\dots,d\}$,
\begin{equation}
\label{eqn:sigma(u_j)}
    \vert\sigma_i(u_j)\vert\asymp_{\eta}\begin{cases}
\Vert \mathbf{Q}\Vert^{j-d-1}& i\leq j-1\\
\Vert\mathbf{Q}\Vert^{j-1}& i\geq j
\end{cases}.
\end{equation}
Denote $\sigma_i(u_j)=u_{ij}$ and let $M_{\mathbf{Q}}=(u_{ij})_{i,j}$ be as in (\ref{eqn:M_Q}). Then, by (\ref{eqn:sigma(u_j)}), (\ref{eqn:M_Q}), and the equality case of the ultrametric inequality,
\begin{equation}
\label{eqn:detM_Q}
\begin{split}
     \vert\det(M_{\mathbf{Q}})\vert=\left|\sum_{\tau\in \mathcal{P}_d}(-1)^{\operatorname{sgn}(\tau)}\prod_{j=1}^d\sigma_j(u_{\tau(j)})\right|=\prod_{j=1}^d\vert \sigma_j(u_j)\vert\\
     =\prod_{i=1}^d\left|\sum_{j=1}^d\sigma_i(u_j)\right|\asymp_{\eta}\Vert \mathbf{Q}\Vert^{\begin{pmatrix}
d\\
2
\end{pmatrix}}.
\end{split}
\end{equation} 
By Theorem \ref{main} and upper-semicontinuity of $\mu$, for every large enough $\Vert \mathbf{Q}\Vert$, $\mu(\mathfrak{x}_{\mathbf{Q}})\leq q^{-d}$. Therefore to prove Theorem \ref{Cassels}, it suffices to show that for every $\Vert \mathbf{Q}\Vert$ large enough, $\mu(\mathfrak{x}_{\mathbf{Q}})\geq q^{-d}$. To do so, we will follow the proof of the main result of \cite{C}. 

Let $\Theta_{\mathbf{Q}}:=\left(\theta_i^{j-1}\right)_{i,j}$ and let 
$$y_{\mathbf{Q}}:=\Theta_{\mathbf{Q}}\left(\mathcal{R}^d+\begin{pmatrix}
\frac{1}{x}\\
\vdots\\
\frac{1}{x}
\end{pmatrix}\right).$$ 
To prove Theorem \ref{Cassels}, it suffices to prove the following theorem.
\begin{theorem}
\label{y_QVal}
For all but finitely many $\mathbf{Q}$, we have that for every $\mathbf{v}\in y_{\mathbf{Q}}$, 
$$N(\mathbf{v})=\prod_{i=1}^d \vert v_i\vert\geq q^{-d}\vert \det(M_{\mathbf{Q}})\vert.$$
\end{theorem}
\begin{proof}[Proof of Theorem \ref{Cassels}]
By definition,
\begin{equation}
\label{eqn:muDef}
    \mu(\mathfrak{x}_{\mathbf{Q}})=\frac{1}{\vert\det(M_{\mathbf{Q}})\vert}\sup_{y\in \operatorname{proj}^{-1}(\mathfrak{x}_{\mathbf{Q}})}\inf_{\mathbf{v}\in y}N(\mathbf{v}).
\end{equation}
By Theorem \ref{y_QVal}, for all but finitely many $\mathbf{Q}$, $\inf_{\mathbf{v}\in y_{\mathbf{Q}}}N(\mathbf{v})\geq q^{-d}\vert \det(M_{\mathbf{Q}})\vert$, and therefore, (\ref{eqn:muDef}) is greater than or equal to $q^{-d}$. On the other hand, by Theorem \ref{main} and upper semicontinuity of $\mu$, $\mu(\mathfrak{x}_{\mathbf{Q}})\leq q^{-d}$, and hence $\mu(\mathfrak{x}_{\mathbf{Q}})=q^{-d}$.
\end{proof}
\begin{lemma}
\label{DiffBase1/x}
$$y_{\mathbf{Q}}=M_{\mathbf{Q}}\left(\mathcal{R}^d+\begin{pmatrix}
\frac{1}{x}\\
\vdots\\
\frac{1}{x}
\end{pmatrix}\right).$$
\end{lemma}
\begin{proof}
Notice that $u_1=1$ and for every $j\in\{2\dots d\}$,
$$u_j=\prod_{l=1}^{j-1}(\theta-Qa_l)=\theta^{j-1}+c_{j,j-2}\theta^{j-2}+\dots +c_{j,1}\theta+c_{j,0},$$
where $c_{j,l}$ are are sums of products of $Qa_1\dots Qa_{j-1}$. Since $a_l\equiv 0\mod x^2$ for every $l\in\{1,\dots ,j-1\}$, then $c_{j,l}\equiv 0 \mod x^2$ for every $j\in\{2,\dots ,d\}$ and for every $l\in\{1,\dots, j-1\}$. Therefore, the change of basis matrix between the $\mathcal{R}$ bases of $\Lambda_{\mathbf{Q}}$,
$\{u_1\dots u_d\}$ and $\{1,\theta,\dots \theta^n\}$, is given by
\begin{equation*}
    P=\begin{pmatrix}
    1&c_{2,0}&\dots&c_{d,0}\\
    0&1&\dots&c_{d,1}\\
    \vdots&\dots&\ddots&\vdots\\
    0&\dots&\dots &1
    \end{pmatrix}.
\end{equation*}
Since $P$ is upper triangular with elements in $x^2\mathcal{R}$ above the diagonal and $c_{i+1,j}\equiv 0 \mod x^2$, then $P^{-1}$ is of the same form as well. Hence, $$P^{-1}\left(\mathcal{R}^d+\begin{pmatrix}
\frac{1}{x}\\
\vdots\\
\frac{1}{x}
\end{pmatrix}\right)\subseteq \mathcal{R}^d+\begin{pmatrix}
\frac{1}{x}\\
\vdots \\
\frac{1}{x}
\end{pmatrix}.$$
Hence,
\begin{equation}
\label{eqn:M_QCont}
\begin{split}
    M_{\mathbf{Q}}\left(\mathcal{R}^d+\begin{pmatrix}
    \frac{1}{x}\\
    \vdots\\
    \frac{1}{x}
    \end{pmatrix}\right)=P\Theta_{\mathbf{Q}}P^{-1}\left(\mathcal{R}^d+\begin{pmatrix}
    \frac{1}{x}\\
    \vdots\\
    \frac{1}{x}
    \end{pmatrix}\right)\subseteq  \Theta_{\mathbf{Q}}\mathcal{R}^d+P\Theta_{\mathbf{Q}}\begin{pmatrix}
    \frac{1}{x}\\
    \vdots \\
    \frac{1}{x}
    \end{pmatrix}.
\end{split}
\end{equation}
Notice that if $\xi=\sum_{i=0}^n\frac{1}{x}\theta^i$, then the fact that $c_{j,l}\equiv 0\mod x^2$ implies that there exist some $b_0\dots b_n\in \mathcal{R}$ such that 
$$P\boldsymbol{\sigma}(\xi)=\sum_{i=0}^d\left(b_i+\frac{1}{x}\right)\boldsymbol{\sigma}(\theta^i)\in \Theta_{\mathbf{Q}}\left(\mathcal{R}^d+\begin{pmatrix}
\frac{1}{x}\\
\vdots\\
\frac{1}{x}
\end{pmatrix}\right).$$ 
Hence, the right hand side of (\ref{eqn:M_QCont}) is contained in $y_{\mathbf{Q}}$. On the other hand, by doing the same procedure and writing $\Theta_{\mathbf{Q}}=P^{-1}M_{\mathbf{Q}}P$, we obtain that $y_{\mathbf{Q}}=M_{\mathbf{Q}}\left(\mathcal{R}^d+\begin{pmatrix}
\frac{1}{x}\\
\vdots\\
\frac{1}{x}
\end{pmatrix}\right)$. 
\end{proof}

Let $\boldsymbol{T}=M_{\mathbf{Q}}^{-1}$. We will need the following bound on the entries of $\boldsymbol{T}$, which can be viewed as an analogue of Lemma 2 in \cite{C}. 
\begin{lemma}
\label{T_ijBnd}
Write $\boldsymbol{T}=(T_{ij})$. Then, $T_{ij}=O(\Vert \mathbf{Q}\Vert^{1-i})$ and $T_{ij}=O(\Vert \mathbf{Q}\Vert^{1-j})$. 
\end{lemma}
\begin{proof}
We use the estimate for the adjugate matrix. Since $\operatorname{adj}(M_{\mathbf{Q}})M_{\mathbf{Q}}=\det(M_{\mathbf{Q}})I$ then,
\begin{equation*}
    \boldsymbol{T}=\operatorname{adj}(M_{\mathbf{Q}})\det(M_{\mathbf{Q}})^{-1}.
\end{equation*}
By the definition of the adjugate matrix $\operatorname{adj}(M_{\mathbf{Q}})_{ij}$ is given by the determinant of the matrix given by removing the $i$-th row and the $j$-th column from (\ref{eqn:M_Q}). Notice that (\ref{eqn:y_QBND}) implies that 
\begin{equation}
\label{eqn:adjj}
\begin{split}
    \left(\operatorname{adj}(M_{\mathbf{Q}})\right)_{ij}\ll \prod_{l\neq j}\vert(M_{\mathbf{Q}})_{ll}\vert=\prod_{l\neq j}\vert \sigma_l(u_l)\vert\\
    \ll \prod_{l\neq j}\Vert \mathbf{Q}\Vert^{l-1}=\Vert \mathbf{Q}\Vert^{\begin{pmatrix}
d\\
2
\end{pmatrix}-(j-1)}.
\end{split}
\end{equation}
Similarly,
\begin{equation}
\label{eqn:adji}
    \left(\operatorname{adj}(M_{\mathbf{Q}})\right)_{ij}\ll \prod_{l\neq i}\Vert \mathbf{Q}\Vert^{l-1}=\Vert \mathbf{Q}\Vert^{\begin{pmatrix}
d\\
2
\end{pmatrix}-(i-1)}.
\end{equation}
Hence, by (\ref{eqn:adjj}),
$$\vert T_{ij}\vert=\vert(M_{\mathbf{Q}}^{-1})_{ij}\vert=\vert \det(M_{\mathbf{Q}})\vert^{-1}\vert \operatorname{adj}(M_{\mathbf{Q}})_{ij}\vert\ll \Vert \mathbf{Q}\Vert^{1-j}.$$
Similarly, by (\ref{eqn:adji}),
$$\vert T_{ij}\vert\ll \Vert \mathbf{Q}\Vert^{1-i}.$$
\end{proof} 
\begin{lemma}
\label{maxXi_i}
There exists some $c>0$ such that for every $\mathbf{v}\in y_{\mathbf{Q}}$, $$\mathbf{v}\notin B(0,c\Vert \mathbf{Q}\Vert^n)\times \dots \times B(0,c\Vert \mathbf{Q}\Vert^n).$$
\end{lemma}
\begin{proof}
Firstly, write $\mathbf{v}=\sum_{i=1}^d\beta_i\boldsymbol{\sigma}(u_i)$ where $\beta_i\equiv \frac{1}{x}\mod \mathcal{R}$. Thus, for every $j\in\{1,\dots,d\}$, $\vert \beta_j\vert \geq \left|\frac{1}{x}\right|=\frac{1}{q}$. Notice that by Lemma \ref{T_ijBnd},
\begin{equation*}
    \frac{1}{q}\leq \vert \beta_d\vert=\left|\sum_{j=1}^d T_{dj}v_j\right|\leq \max_j\vert v_j\vert\cdot \vert T_{dj}\vert\ll \Vert \mathbf{Q}\Vert^{-n}\max_j\vert v_j\vert.
\end{equation*}
Thus, there exists some $c>0$ such that $\max_j\vert v_j\vert>c\Vert \mathbf{Q}\Vert^{n}$ and hence the claim follows. 
\end{proof}
\begin{lemma}
\label{staby_Q}
For every $i\in\{1\dots d\}$, $\operatorname{diag}(\boldsymbol{\sigma}(\omega_i))\in \operatorname{stab}_{A_1}(y_{\mathbf{Q}})$.
\end{lemma}
\begin{proof}
We shall now show that for every $j\in\{1,\dots ,d\}$, the unit $\omega_j$ preserves the grid 
$$\iota_{\mathbf{Q}}:=\Bigg\{\sum_{i=1}^d\left(b_i+\frac{1}{x}\right)\theta^{i-1}:b_i\in \mathcal{R}\Bigg\}.$$
This will imply that $\operatorname{diag}\left(\boldsymbol{\sigma}(\omega_j)\right)$ preserves the grid $y_{\mathbf{Q}}=\boldsymbol{\sigma}(\iota_{\mathbf{Q}})$. If $\xi=\sum_{i=0}^n\alpha_i\theta^i\in \iota_{\mathbf{Q}}$ then, for every $j\in\{1,\dots, d\}$, 
\begin{equation*}
    (\theta-Qa_j)\sum_{i=0}^n\alpha_i\theta^i=\alpha_n\theta^d+\sum_{i=1}^{n}\theta^i(\alpha_{i-1}-Qa_j\alpha_i)-Qa_j\alpha_0.
\end{equation*}
Due to (\ref{eqn:P_QDef}), $\theta^d=c_0+c_1\theta+\dots +c_n\theta^n$ where $c_i$ are composed of sums of products of the polynomials $Qa_i$ and $c_0=1+Q^da_1\dots a_d$. In particular $c_i\equiv 0\mod x^2$ for every $i\in\{1\dots n\}$ and $c_0\equiv 1 \mod x^2$. Hence, 
\begin{equation*}
    (\theta-Qa_j)\sum_{i=0}^n\alpha_i\theta^i=\sum_{i=1}^n\theta^i(\alpha_{i-1}-Qa_j\alpha_i+c_i\alpha_n)+(c_0\alpha_n-Qa_j\alpha_0).
\end{equation*}
Since $\alpha_i\equiv \frac{1}{x}\mod \mathcal{R}$, $a_j\equiv 0\mod x^2$, and $c_i\equiv 0\mod x^2$ then, $Qa_j\alpha_i\equiv 0\mod \mathcal{R}$ and $c_i\alpha_n\equiv 0\mod \mathcal{R}$ for every $i\in\{1,\dots ,n\}$. Therefore, for every $i\in\{1\dots n\}$, 
$$\alpha_{i-1}-Qa_j\alpha_i+c_i\alpha_n\equiv \alpha_{i-1}\mod \mathcal{R}\equiv \frac{1}{x}\mod \mathcal{R}.$$
In addition, since $c_0=1+Q^da_1\dots a_d$, then, 
$$c_0\alpha_n-Qa_j\alpha_0\equiv \alpha_n\mod \mathcal{R}\equiv \frac{1}{x}\mod \mathcal{R}.$$
Hence, for every $j\in\{1,\dots, d\}$, the unit $\omega_j$ preserves the grid $\iota_{\mathbf{Q}}$, so that $\operatorname{diag}\left(\boldsymbol{\sigma}(\omega_j)\right)\in \operatorname{stab}_{A_1}(y_{\mathbf{Q}})$.
\end{proof}
To conclude the proof of Theorem \ref{Cassels}, we shall use the following proposition.
\begin{proposition}
\label{xiWorks}
For every $0<\varepsilon<\frac{1}{2}$, there exists some $C'>0$ such that for all but finitely many $\mathbf{Q}$, and for every $\mathbf{v}\in y_{\mathbf{Q}}$ such that
\begin{equation}
\label{eqn:prodsmall}
    \prod_{i=1}^d\vert v_i\vert\leq q^{-d}\vert \det(M_{\mathbf{Q}})\vert,
\end{equation}
there exist some $b_1\dots b_n\in \mathbb{Z}$ and a permutation $\tau\in \mathcal{P}_d$ such that the vector 
$$\mathbf{v}'=\operatorname{diag}(\boldsymbol{\sigma}(\omega_1))^{b_1}\dots \operatorname{diag}(\boldsymbol{\sigma}(\omega_n))^{b_n}\mathbf{v}$$
satisfies
    \begin{equation}
       \label{eqn:vinBox} \mathbf{v}'\in B\left(0,C'\Vert \mathbf{Q}\Vert^{\tau(1)-\varepsilon}\right)\times B\left(0,C'\Vert \mathbf{Q}\Vert^{\tau(2)-\varepsilon}\right)\times \dots \times B\left(0,C'\Vert \mathbf{Q}\Vert^{\tau(d)-\varepsilon}\right).
    \end{equation}
\end{proposition}
In addition, we shall use the following lemma. 
\begin{lemma}
\label{ProdLarge}
There exists some $C'>0$ such that for all but finitely many $\mathbf{Q}$, if $\mathbf{v}$ satisfies (\ref{eqn:vinBox}) for some $\tau\in \mathcal{P}_d$ and some $\varepsilon\in (0,1)$ then, 
\begin{equation}
\label{eqn:ProdXi}
    \prod_{i=1}^d\vert v_i\vert\geq q^{-d}\vert \det(M_{\mathbf{Q}})\vert.
\end{equation}
\end{lemma}
\begin{proof}
Without loss of generality, assume that $\tau=Id$ and write $\mathbf{v}=\sum_{i=1}^d\beta_i\boldsymbol{\sigma}(u_i)$. By Lemma \ref{T_ijBnd} and (\ref{eqn:vinBox}),
\begin{equation}
\label{eqn:beta_iBnd}
\begin{split}
    \vert \beta_i\vert=\left|\sum_{j=1}^dT_{ij}v_j\right|\leq \max_{j=1,\dots d} \vert T_{ij}\vert \cdot \vert v_j\vert\\
    \ll \max_{j=1,\dots d} O(\Vert \mathbf{Q}\Vert^{1-j+j-\varepsilon})=O(\Vert \mathbf{Q}\Vert^{\varepsilon}).
\end{split}
\end{equation}
Thus, by (\ref{eqn:sigma(u_j)}) and (\ref{eqn:beta_iBnd}),
\begin{equation}
\label{eqn:xi_i-beta_iu_ii}
\begin{aligned}
    \vert v_i-\beta_i\sigma_i(u_i)\vert=\left|\sum_{j\neq i}\beta_j\sigma_i(u_j)\right|\\
    \leq \max_{j\neq i}\vert\beta_j\vert \cdot \vert \sigma_i(u_j)\vert\ll \Vert \mathbf{Q}\Vert^{\varepsilon}\Vert \mathbf{Q}\Vert^{i-2}=O(\Vert \mathbf{Q}\Vert^{i-1-\varepsilon}).
\end{aligned}
\end{equation}
Since $\mathbf{v}\in y_{\mathbf{Q}}$, them $\beta_i\equiv\frac{1}{x}
\mod \mathcal{R}$. Hence, for every $i\in\{1,\dots, d\}$, we have that $\vert \beta_i\vert\geq \frac{1}{q}$. Thus, by the equality case of the ultrametric inequality (\ref{eqn:ultrametric}), (\ref{eqn:sigma(u_j)}), (\ref{eqn:beta_iBnd}), and (\ref{eqn:xi_i-beta_iu_ii})  we obtain that 
\begin{equation}
\label{eqn:xi_iBnd}
    \vert v_i\vert=\max\{\vert v_i-\beta_i\sigma_i(u_i)\vert,\vert \beta_i\sigma_i(u_i)\vert\}\geq \max\Bigg\{O(\Vert \mathbf{Q}\Vert^{i-1-\varepsilon}), \frac{1}{q}\vert \sigma_i(u_i)\vert\Bigg\}
\end{equation}
Hence, (\ref{eqn:xi_iBnd}) and (\ref{eqn:detM_Q}) imply that
\begin{equation*}
    \prod_{i=1}^d\vert v_i\vert\geq (1+o(1))\prod_{i=1}^d\frac{1}{q}\vert \sigma_i(u_i)\vert=\frac{1}{q^d}\vert \det(M_{\mathbf{Q}})\vert.
\end{equation*}
\end{proof}
\begin{proof}[Proof of Theorem \ref{y_QVal}]
Let $\mathbf{v}\in y_{\mathbf{Q}}$. If 
\begin{equation*}
    \prod_{i=1}^d\vert v_i\vert>\frac{1}{q^d}\vert \det(M_{\mathbf{Q}})\vert,
\end{equation*}
then there is nothing to check. Hence, we shall assume that $\mathbf{v}$ satisfies (\ref{eqn:prodsmall}). By Proposition \ref{xiWorks}, there exist some $b_1\dots b_n\in \mathbb{Z}$ such that $$\mathbf{v}'=\operatorname{diag}(\boldsymbol{\sigma}(\omega_1))^{b_1}\dots \operatorname{diag}(\boldsymbol{\sigma}(\omega_n))^{b_n}\mathbf{v}$$
satisfies (\ref{eqn:vinBox}) for some $\tau\in \mathcal{P}_d$. Hence, Lemma \ref{staby_Q} and Lemma \ref{ProdLarge} imply together that 
$$N(\mathbf{v})=N(\mathbf{v}')\geq q^{-d}\vert \det(M_{\mathbf{Q}})\vert.$$
Thus, by discreteness of the value set of $N$ around non-zero points, we obtain that for all but finitely many $\mathbf{Q}$, every $\mathbf{v}\in y_{\mathbf{Q}}$ satisfies $N(\mathbf{v})\geq q^{-d}\vert \det(M_{\mathbf{Q}})\vert$.
\end{proof}
Hence to conclude the proof of Theorem \ref{y_QVal}, it suffices to prove Proposition \ref{xiWorks}. 
\subsection{Proof of Proposition \ref{xiWorks}}
Let $\mathbf{v}\in y_{\mathbf{Q}}$ be such that $\mathbf{v}$ satisfies (\ref{eqn:prodsmall}). Then, there exist $t\in \{0,1\dots n\}$ and $s\in \mathbb{Z}$ such that $\vert v_1\dots v_d\vert=q^t\cdot q^{ds}$. Thus, there exists a diagonal matrix $\mathbf{g}=\operatorname{diag}\{x^t,1\dots 1\}$ and a vector $\mathbf{v}^{(0)}$ which satisfies $\prod_{i=1}^d\left|\mathbf{v}^{(0)}_i\right|=1$, such that
\begin{equation}
\label{eqn:XiMatrix}
    \mathbf{v}=x^s\mathbf{g}\mathbf{v}^{(0)}=x^s\begin{pmatrix}
    x^t&&&\\
    &1&&\\
    &&\ddots&\\
    &&&1
    \end{pmatrix}\mathbf{v}^{(0)}.
\end{equation}
Since $\mathbf{v}$ satisfies (\ref{eqn:prodsmall}),
\begin{equation}
\label{eqn:smallProd}
    q^{sd+t}=N(\mathbf{v})=\prod_{i=1}^d\vert v_i\vert\leq   \frac{1}{q^d}\vert \det(M_{\mathbf{Q}})\vert.
\end{equation}

\begin{theorem}
\label{rin[0,1]}
Let $b_1,\dots b_d\in \mathbb{Z}$ and define $$\tilde{\mathbf{v}}^{(0)}=\left(\prod_{j=1}^d\operatorname{diag}(\boldsymbol{\sigma}(\omega_j))^{b_j}\right)\mathbf{v}^{(0)}.$$
Then,
\begin{enumerate}
    \item $\prod_{i=1}^d\left| \tilde{v}^{(0)}_i\right|=\prod_{i=1}^d\left| v^{(0)}_i\right|=1$. 
    \item \label{viaV'} Moreover, there exists some vector $\tilde{\mathbf{v}}\in y_{\mathbf{Q}}$, such that the vector $\tilde{\mathbf{v}}^{(0)}$ arises from $\tilde{\mathbf{v}}$ through (\ref{eqn:XiMatrix}).
\end{enumerate}
\end{theorem}
\begin{proof}
Since $\omega_j$ are units in $\mathcal{O}_{\mathbb{F}_{\mathbf{Q}}}$, then $\operatorname{diag}(\boldsymbol{\sigma}(\omega_j))\in A_1$ for every $j\in\{1,\dots,d\}$. Hence,
\begin{equation*}
    \prod_{i=1}^d\left|\tilde{v}^{(0)}_i\right|=\prod_{i=1}^d\left| v^{(0)}_i\right|\cdot \prod_{j=1}^d\vert\theta_i-Qa_j\vert^{b_j}=\prod_{i=1}^d\left| v^{(0)}_i\right|=1.
\end{equation*}
Define $\tilde{\mathbf{v}}=\operatorname{diag}\left(\boldsymbol{\sigma}\left(\prod_{j=1}^d\omega_j^{b_j}\right)\right)\mathbf{v}$. Then, by Lemma \ref{staby_Q}, $\tilde{\mathbf{v}}\in y_{\mathbf{Q}}$. Moreover,
\begin{equation*}
\begin{split}
    \tilde{\mathbf{v}}=\operatorname{diag}\left(\prod_{j=1}^d\sigma(\omega_j)^{b_j}\right)\mathbf{v}
    =x^s\left(\prod_{j=1}^d\sigma(\omega_j)^{b_j}\right)\mathbf{g}\mathbf{v}^{(0)}=x^s\mathbf{g}\tilde{\mathbf{v}}^{(0)}.
\end{split}
\end{equation*}
\end{proof}
We shall now reinterpret Proposition \ref{CoveringR} for real simplex sets very close to the standard simplex set. Let $\Phi_*$ be the standard simplex set (see (\ref{eqn:stdSimp})) and let $\psi_*=\rho(\Phi_*)\subseteq \mathbb{R}_0^d$. Then we can reinterpret Proposition \ref{CoveringR} or equivalently the corollary to Lemma 2 in \cite{C} in the following manner. 
\begin{lemma}
\label{covAlmStd}
For any $\varepsilon>0$, there exists $\delta>0$ such that if $\psi\subseteq \mathbb{R}_0^d$ is a simplex set satisfying $\psi=T\psi_*$, where $T\in \operatorname{GL}_d(\mathbb{R})$ satisfies $\Vert T-I\Vert<\delta$ then,
\begin{enumerate}
    \item \label{eqn:covAlmStdAll} $\langle \psi\rangle+\Bigg\{\mathbf{u}\in \mathbb{R}_0^d:\lceil \mathbf{u}\rceil_{\mathbb{R}_0^d}\leq \frac{n}{2}(1+\varepsilon)\Bigg\}=\mathbb{R}_0^d$ and
    \item \label{eqn:covAlmStd} $\mathbb{R}_0^d\setminus \left(\langle \psi\rangle+\bigg\{\mathbf{u}\in \mathbb{R}_0^d:\lceil \mathbf{u}\rceil_{\mathbb{R}_0^d} \leq \frac{n}{2}(1-\delta)\Bigg\}\right)\subseteq B_{\varepsilon}(W_*)+\langle \psi\rangle$,
\end{enumerate}
where $W_*$ is a set of $n!$ vectors obtained by permuting the coordinates of 
\begin{equation*}
    \mathbf{w}=\frac{1}{d}\sum_{l=1}^n(l-1)\rho(\mathbf{t}_l)=\begin{pmatrix}
    \frac{n}{2}\\
    \frac{n}{2}-1\\
    \vdots\\
    -\frac{n}{2}
    \end{pmatrix}.
\end{equation*}
\end{lemma}
\begin{remark}
Lemma \ref{covAlmStd} holds, since as $\psi\rightarrow \psi_*$, $W_{\psi}\rightarrow W_*$ and $S_{\psi}\rightarrow S_*$. 
\end{remark}
\begin{proof}[Proof of Proposition \ref{xiWorks}]
Let $\varepsilon<\frac{1}{2}$. Let $\mathbf{v}\in y_{\mathbf{Q}}$ satisfy (\ref{eqn:smallProd}). Then by (\ref{eqn:XiMatrix}), there exist $s\in \mathbb{Z}$, a diagonal matrix $\mathbf{g}=\operatorname{diag}\{x^t,1,\dots,1\}$, where $0\leq t\leq n$, and $\mathbf{v}^{(0)}$ with $\prod_{i=1}^d\left|v^{(0)}_i\right|=1$ such that $\mathbf{v}=x^s\mathbf{g}\mathbf{v}^{(0)}$. Furthermore, by Lemma \ref{maxXi_i},
\begin{equation}
\label{eqn:xi0ProdLarge}
    \Vert\mathbf{Q}\Vert^n\ll \max_{i=1,\dots,d}\vert v_i\vert=q^s\max_{i=1,\dots ,d}\vert \mathbf{g}_{ii}\vert\cdot \left|v_i^{(0)}\right|\leq q^{s+n}\max_{i=1\dots d}\left|v_i^{(0)}\right|.
\end{equation}
By (\ref{eqn:detM_Q}) and (\ref{eqn:smallProd}),
\begin{equation}
\label{eqn:q^-sBnd}
    q^{-s}\geq q^{\frac{t}{d}+1}\vert \det(M_{\mathbf{Q}})\vert^{-\frac{1}{d}}\geq q\vert \det(M_{\mathbf{Q}})\vert^{-\frac{1}{d}}\gg \Vert \mathbf{Q}\Vert^{-\frac{n}{2}}.
\end{equation}
Thus, (\ref{eqn:xi0ProdLarge}) and (\ref{eqn:q^-sBnd}) imply that
\begin{equation}
\label{eqn:maxxii}
     \max_{i=1,\dots d}\left|v^{(0)}_i\right|\gg \Vert \mathbf{Q}\Vert^{\frac{n}{2}}.
\end{equation}
By (\ref{eqn:samej}) and (\ref{eqn:diffj}), the simplex set $$\psi_{\mathbf{Q}}=\Bigg\{\frac{1}{\log\Vert \mathbf{Q}\Vert}\rho\left(\boldsymbol{\sigma}(\omega_1)\right)\dots \frac{1}{\log\Vert \mathbf{Q}\Vert}\rho\left(\boldsymbol{\sigma}(\omega_d)\right)\Bigg\}$$
converges to the simplex set $\psi_*$. Thus, Lemma \ref{covAlmStd}(\ref{eqn:covAlmStdAll}) implies that whenever $\Vert\mathbf{Q}\Vert$ is large enough, there exists some unit $\omega\in \mathcal{O}_{\mathbb{F}_\mathbf{Q}}$ such that $$\Bigg\lceil\frac{1}{\log\Vert \mathbf{Q}\Vert}\rho\left(\mathbf{v}^{(0)}\right)+\frac{1}{\log\Vert \mathbf{Q}\Vert}\rho(\boldsymbol{\sigma}(\omega))\Bigg\rceil\leq \frac{n}{2}+\varepsilon.$$
Since $\operatorname{\sigma}(\omega)$ preserves $y_{\mathbf{Q}}$, we can assume that $\omega=1$. On the other hand, (\ref{eqn:maxxii}) implies that whenever $\Vert \mathbf{Q}\Vert$ is large enough, $\frac{1}{\log\Vert \mathbf{Q}\Vert}\rho\left(\mathbf{v}^{(0)}\right)$ belongs to the left hand side of Lemma \ref{covAlmStd}(\ref{eqn:covAlmStd}), where $\delta$ is the $\delta$ corresponding to $\varepsilon$ from Lemma \ref{covAlmStd}. Hence, Lemma \ref{covAlmStd}(\ref{eqn:covAlmStd}) implies that there exists a units vector $\mathbf{b}=\rho\left(\prod_{i=1}^d\boldsymbol{\sigma}(\omega_i)^{d_i}\right)$ where $d_i\in \mathbb{Z}$ and a permutation $\tau\in \mathcal{P}_d$, such that 
\begin{equation}
\label{eqn:vcloseW}
    \Bigg\lceil\frac{1}{\log\Vert \mathbf{Q}\Vert}\rho\left(\mathbf{v}^{(0)}\right)+\frac{1}{\log\Vert \mathbf{Q}\Vert}\rho(\mathbf{b})-\tau(\mathbf{w})\Bigg\rceil<\varepsilon.
\end{equation}
Since $\boldsymbol{\sigma}(\omega_i)$ preserve $y_{\mathbf{Q}}$ for every $i\in\{1,\dots,d\}$, we can assume for simplicity that $\mathbf{b}=Id$ and that $\tau=Id$. Then, (\ref{eqn:vcloseW}) implies that
\begin{equation*}
    \Bigg\lceil \rho\left(\mathbf{v}^{(0)}\right)-\begin{pmatrix}
    \frac{n}{2}\\
    \frac{n}{2}-1\\
    \vdots \\
    -\frac{n}{2}
    \end{pmatrix}\log\Vert \mathbf{Q}\Vert \Bigg\rceil<\varepsilon\log\Vert \mathbf{Q}\Vert.
\end{equation*}
Hence, for every $i\in\{1\dots d\}$,
\begin{equation*}
    \rho\left(v_i^{(0)}\right)\leq \left(\frac{n}{2}-i+\varepsilon\right)\log\Vert \mathbf{Q}\Vert.
\end{equation*}
Therefore, 
\begin{equation}
\label{eqn:v^0_i}
    \left| v_i^{(0)}\right|\leq\Vert \mathbf{Q}\Vert^{\frac{n}{2}-i+\varepsilon}.
\end{equation}
Thus, (\ref{eqn:q^-sBnd}) and (\ref{eqn:v^0_i}) imply that
\begin{equation}
\label{eqn:viBND}
    \vert v_i\vert\leq q^{s+n}\left|v_i\right|\ll \Vert \mathbf{Q}\Vert^{n-i+\varepsilon}.
\end{equation}
Define $\tau'(i):=\begin{cases}
d-i&i=1\dots n\\
d&i=d
\end{cases}$. Then since $\varepsilon\leq \frac{1}{2}$, $n-i+\varepsilon\leq \tau'(i)-\varepsilon$. Thus, by (\ref{eqn:viBND}), 
\begin{equation*}
    \vert v_i\vert\ll \Vert \mathbf{Q}\Vert^{\tau'(i)-\varepsilon}.
\end{equation*}
This shows that for every large enough $\Vert \mathbf{Q}\Vert$ and for every $\mathbf{v}$ satisfying (\ref{eqn:smallProd}), there exists a diagonal matrix $\prod_{i=1}^d\boldsymbol{\sigma}(\omega_i)^{m_i}$ such that $\mathbf{u}=\prod_{i=1}^d\boldsymbol{\sigma}(\omega_i)^{m_i}\mathbf{v}$ satisfies (\ref{eqn:vinBox}). This concludes the proof of Proposition \ref{xiWorks} and hence the proof of Theorem \ref{Cassels}. 
\end{proof}
\bibliography{Ref}
\bibliographystyle{amsalpha}
\end{document}